\documentclass[leqno,11pt,letterpaper, english]{amsart}
\usepackage[usenames,dvipsnames]{color}
\usepackage[colorlinks=true,linkcolor=Red,citecolor=Green]{hyperref}
\usepackage{amsmath,amssymb,amsthm,graphicx,url}
\usepackage{epstopdf}
\usepackage{color}
\usepackage{dsfont}
\usepackage[T1]{fontenc}
\usepackage{natbib}
\usepackage{comment}

\newcommand{\xqedhere}[1]{%
	\rlap{%
		\hbox to#1{%
			\hfil
			\llap{%
				\ensuremath{\square}
			}%
		}%
	}%
}

\def\pasdegrille{\let\grille = \pasgrille}

\def\aat#1#2#3{
	\divide \dimen1 by 48 \dimen3=\dimen1 \multiply \dimen1 by #1
	\advance \dimen1 by -\dimen3 \divide \dimen1 by 101 \multiply
	\dimen1 by 100 \divide \dimen2 by \count11 \multiply \dimen2 by #2
	\setbox0=\hbox{#3}\ht0=0pt\dp0=0pt
	\rlap{\kern\dimen1 \vbox to0pt{\kern-\dimen2\box0\vss}}\dimen1= \wd1
	\dimen2=\ht1}
\def\pasgrille{
	\count12= \dimen1 \divide \count12 by 50 \divide \dimen2 by
	\count12 \count11 =\dimen2 \ \divide \dimen1 by 48
	\setlength{\unitlength}{\dimen1} \smash{\rlap{\ }} \dimen1= \wd1
	\dimen2=\ht1 }
\def\grille{
	\count12= \dimen1 \divide \count12 by 50 \divide \dimen2 by
	\count12 \count11 =\dimen2 \ \divide \dimen1 by 48
	\setlength{\unitlength}{\dimen1}
	\smash{\rlap{\graphpaper[1](0,0)(50, \count11)}} \dimen1= \wd1
	\dimen2=\ht1 }

\pasdegrille
\newcommand{\be}{\begin{equation}}
	\newcommand{\ee}{\end{equation}}

\renewcommand{\Re}{\mathop{\rm Re}\nolimits}

\newcommand{\nn}{\nonumber} 

\theoremstyle{plain}

\newtheorem{thm}{Theorem}

\newtheorem{prop}{Proposition}[section]
\newtheorem{cor}[prop]{Corollary}
\newtheorem{lem}[prop]{Lemma}

\theoremstyle{definition}

\numberwithin{equation}{section}

\def\squarebox#1{\hbox to #1{\hfill\vbox to #1{\vfill}}}

\usepackage{amsxtra}

\ifx\pdfoutput\undefined
\DeclareGraphicsExtensions{.pstex, .eps}
\else
\ifx\pdfoutput\relax
\DeclareGraphicsExtensions{.pstex, .eps}
\else
\ifnum\pdfoutput>0
\DeclareGraphicsExtensions{.pdf}
\else
\DeclareGraphicsExtensions{.pstex, .eps}
\fi
\fi
\fi

\title[Unique continuation for semiclassical Witten-Schrödinger type operator ]{Quantitative unique continuation for semiclassical Witten-Schrödinger type operator with low regular coefficients.}

\author[S.~Campagne]{Sébastien Campagne}
\address{Universit{\'e} Paris-Saclay, Math{\'e}matiques, UMR 8628 du CNRS, B{\^a}t 307, 91405  Orsay Cedex, France,   and Institut Universitaire de France}
\email{sebastien.campagne@universite-paris-saclay.fr}

\date{\today}

\usepackage{amssymb}
\usepackage{amsmath, amsthm, amsopn, amsfonts}

\def\11{{\rm 1~\hspace{-1.4ex}l} }
\def\R{\mathbb R}
\def\C{\mathbb C}

\def\N{\mathbb N}

\begin{document}   
	
	\begin{abstract}
		In this note we study eigenfunctions localization of the semiclassical selfadjoint Schrödinger operator with low regular coefficients on domains of $\R^2$. We show a unique continuation result on the local decay of the eigenfunction: local supremum on open sets can't vanish faster than global supremum on these same sets times an exponential term that decrease in the semiclassical limit $h\to 0$ as $\exp(-c\log(1/h)/h)$ in the worst case. We deduce from this result  a control on the eigenfunctions of the semiclassical Schrödinger operator on a $W^{2,\infty}$ compact Riemann surface with Lipschitz positive metric and with Robin boundary condition. This last part leads us to study the Schrödinger operator with a Witten Laplacian.

		\begin{center} \rule{0.6\textwidth}{.4pt} \end{center}\end{abstract}     
	
	\ \vskip -1cm \noindent\hfil\rule{0.9\textwidth}{.4pt}\hfil \vskip 1cm 
	
	\maketitle
	
	

	\section{Introduction}
	
	In this article, we study the localization properties of eigenfunctions $u_h$ for  semiclassical Witten-Schrödinger type operator $P_h$ on $\R^2$. Specifically, we consider the eigenvalue equation:
	
	\begin{equation}
		\label{eq: Schrodinger equation}
		P_hu_h:= \left(-\frac{h^2}{K_h}\mathrm{div}(K_h\nabla) +V_h\right)u_h=E_hu_h \quad \mbox{ on } B(0,R)\subset \R^2,
	\end{equation}
	
	where 
	\begin{itemize}
		\item $K_h:=\exp(\phi_h/h) \mbox{ with } \phi_h=\phi+\mathcal{O}_{W^{1,\infty}}(h) \in W^{1,\infty}(B(0,R),\R),$
		\item $V_h:=V+\mathcal{O}_{L^\infty}(h)\in L^\infty(B(0,R),\R)$,
		\item and $E_h:=E+\mathcal{O}_{\R}(h)\in \R$,
	\end{itemize}
	 with $h$ the semiclassical parameter. If the function $\phi_h$ is sufficiently regular, we can rewrite the operator $P_h$ in the form $-h^2\Delta_{\phi_h/2}+\tilde{V}$, where $\Delta_{\phi_h/2}$ is the Witten Laplacian $e^{-\phi_h/2h}\Delta e^{\phi_h/2h}$. This deformation of the usual Laplacian was introduce by \cite{10.4310/jdg/1214437492} to study the Morse theory and Supersymmetry. It also appeared in some study of quantum tunneling  and metastable states as for example in the Kramers-Smoluchowski equation (see for instance \cite*{Kramers-Smoluchowski}).

 Our goal here is to understand how the weight $K_h$ and the potential $V_h$ affect the spatial concentration of $u_h$ in the semiclassical limit $h \rightarrow 0$.

	Our first result in this article is the following semiclassical adaptation of Theorem $2.2$ from \cite{logunov2020landisconjectureexponentialdecay} generalised to pondered Schrödinger operator:
	\begin{thm}
		\label{thm: mainthm}
		For $h$ small enough, let $u_h\in H^1(B(0,R),\R)$ be a solution on  $B(0,R)$ of the equation $$P_{h}u_h:=\left(-\frac{h^2}{K_h}\mathrm{div}(K_h\nabla)+ V_h\right)u_h=E_h u_h $$ in the sense of distributions:
		for all $\chi \in \mathcal{C}_0^\infty(B(0,R),\R)$, $$\int \left(h^2\nabla u_h\cdot\nabla \chi +V_hu_h\chi\right)K_h\mathrm{d}x=\int E_hu_h\chi K_h\mathrm{d}x.$$ 
		Let $R_0<R$ and suppose that 
		$$\sup_{B(0,R_0/2)} |u_h|\geq e^{-N_h}  \sup_{B(0,R_0)} |u_h|,$$
		then for $0<r\leq R_0/4$:
		\begin{equation}
			\sup \limits _{B(0,r)} |u_h| \geq \left(\frac{r}{R_0} \right)^{ L\left(\frac{R_0}{h}\epsilon(h) + N_h+\frac{R_0}{h}	\mathrm{osc}(\phi_h)\right)}\,\sup \limits_{B(0,R_0)}|u_h|,
		\end{equation}
		with $L$ an absolute constant, $\mathrm{osc}(\phi_h)= \max \phi_h -\min \phi_h$ and
		$$
		\epsilon(h):=  \left\{
		\begin{array}{ll}
			\Theta\log\left(\frac{R_0\Theta}{h}\right) & \mbox{if } 	||\nabla \phi||_{\infty}>0,\\
			\Theta\log\left(\frac{R_0\Theta}{h}\right)^{1/2} & \mbox{else.}
		\end{array}
		\right.
		$$ 
		where $\Theta= \max \{\|\nabla \phi\|_{\infty}, \sqrt{\|V-E\|_{\infty}},1\}$. 
	\end{thm}
	
	If we forget the exponentials weights, we recover the usual semiclassical Schrödinger operator on $\R^2$. Thus this result gives constraints on speed limit of the Schrödinger eigenfunctions' localization in a reasonable neighbourhood of the maximum. This works as a counterpart of the usual Agmon estimate that gives at least exponential decrease inside the classical forbidden region $\{V> E\}$ (see for instance \cite*{Helffer1988}, Chapter $7$). This type of unique continuation is generally achieved thanks to Carleman inequality for smooth elliptic operator. In this case, Carleman inequality shows that eigenfunctions can't localize faster than exponentially.  With the low regularity of our parameters, Carleman inequality can't be used any more. In fact, it still works in dimension $1$ as shown by \cite*{dyatlovMathematicalTheoryScattering2019} (see Theorem $2.32$). In the same dimension \cite*{laurentUniformObservationSemiclassicalSchrodinger2023} even show a unique continuation property for some class of potential in term of Agmon distance (see Theorem $1.4$). In higher dimension, 
	\cite*{kloppSemiclassicalResolventEstimates2019a} constructed a Carleman inequality  for bounded complex potential (lemma $4$). Their result imply a maximal rate of localization of $e^{-Ch^{-4/3}}$, which is faster than the usual Carleman one. If we add some regularity to the potential, we can find better estimates as done by  \cite*{vodevSemiclassicalResolventEstimates2020}  who obtain a Carleman estimate for $\alpha-$Hölder regular potentials in term of $\alpha$  ($e^{-Ch^{-4/(\alpha+3)}}$ , Theorem $3.1$). We can also assume others properties for the potential as for instance \cite*{datchev2014quantitativelimitingabsorptionprinciple} in dimension $\geq 3$ who ask a radial decrease control on the potential and its radial derivative (Lemma $2.2$) or \cite{shapiro2017semiclassicalresolventboundsdimension} in dimension $2$ who ask a decrease control on the gradient of the potential (Lemma $3.2$).
	
	The proof of the Theorem \ref{thm: mainthm} use the method developed by \cite*{logunov2020landisconjectureexponentialdecay} to study the Landis conjecture in dimension $2$.
	As a reminder, the Landis conjecture (see \cite{Kondratev1991}) deals with how quickly solutions $u \neq 0$ to
	
	\begin{equation}
		(-\Delta+V)u=0
	\end{equation}
	can decay at infinity when $V$ is bounded.
	Landis and Kondrat'ev conjectured that non-zero solutions cannot decay faster than exponentially. This was disproven by \cite{Meshkov1992} for complex potentials: he found a bounded potential and a solution that decay as $e^{-C|x|^{4/3}}$ in dimension $2$. He also proved that is was impossible for a non zero eigenfunction to decrease faster than $e^{-|x|^{\epsilon+4/3}}$ in any dimension higher than $1$. Later, \cite*{Bourgain2005} obtained a quantitative estimate of the Meshkov result  (Lemma $3.10$). They showed that Schrödinger eigenfunctions on $\R^d$ could not decrease faster than $ e^{-C(|x|^{4/3}\log|x|} .$ 
	Naturally, the question of whether the conjecture was true in the case of a real potential was raised. Some results in this regard have been obtained. First in dimension $2$, \cite*{Davey2017} proved a weaker version of the Landis conjecture for second order elliptic operator with variable coefficients and with potential $V$ assumed to be non-negative.  They obtained a decrease speed limit of $e^{-C|x|\log(|x|)}$.
	Their result was improve by \cite*{logunov2020landisconjectureexponentialdecay} in case of the usual Schrödinger operator (not assumed to be non-negative), where they obtained a maximal decay of $e^{-CR\log(R)^{1/2}}$.
	Their proof involves removing holes in the domain to increase the control of how functions oscillate by reducing the Poincaré constant. Then they use quasi-conformal mappings to transform the problem into one about harmonic functions on a domain with holes.
	To go deeper, \cite*{fernandez-bertolinLandisConjectureSurvey2024} did a complete survey of Landis' conjecture. 
	Unfortunately, the conjecture proved to be false in dimension $\geq 3$. In recent work of \cite{frank2026counterexampleslandisconjecturedimensions}, a counterexemple 
	decaying at speed $e^{-C|x|^{4/3}}$ was constructed thank to AI, thus recovering the result of Meskov. However, it seems that the conjecture is true if we consider non-negative operator as shown by \cite*{das2024landisconjectureliouvillecomparison}. Their result even show that the maximum rate of decline is slightly slower than Landis had predicted. 
	
	In a second time, we use the result of Theorem \ref{thm: mainthm} to get information on the localization of Schrödinger's eigenfunctions in a compact Riemannian setting. More precisely, let $(M,g)$ be a $W^{2,\infty}$ Riemann surface ( $W^{2,\infty}$ means that the changing of charts are $W^{2,\infty}\cap \mathcal{C}^1$) equipped with $g$ a Lipschitz positive metric. We deduce from Theorem \ref{thm: mainthm} the following result on closed (compact and without boundary) surface: 
	
	\begin{thm}
		Assume that $M$ is a closed Riemannian surface. For $h$ small enough, let be  $u_h\in H^1(M,\R)$ weak solution of $$-\frac{h^2}{K_h}\mathrm{div}_g( K_h\nabla_g )u_h + V_h u_h =E_hu_h \mbox{ on } M,$$ where  $K_h=e^{\phi_h/h}$ with $\phi_h=\phi+\mathcal{O}_{W^{1,\infty}}(h)\in W^{1,\infty}(M,\R)$, $V_h=V+\mathcal{O}_{L^\infty}(h)\in L^\infty(M,\R)$ and $E_h= E +\mathcal{O}_{\R}(h) \in\R$. 
		
		Then for any  ball $B_r$ of radius $r<1/4$:
		
		\begin{equation}
			\sup\limits_{B_r}|u_h| \geq  r^{ L\left(\frac{\epsilon(h)+\mathrm{osc}(\phi_h)}{h}\right)} \sup_{M}|u_h|,
		\end{equation}
		where $L$ is an absolute constant, $\mathrm{osc}(\phi_h)= \max \phi_h -\min \phi_h$  and 
		$$
		\epsilon(h):=  \left\{
		\begin{array}{ll}
			\Theta\log\left(\frac{\Theta}{h}\right) & \mbox{if } 	||\nabla \phi||_{\infty}>0,\\
			\Theta\log\left(\frac{\Theta}{h}\right)^{1/2} & \mbox{else.}
		\end{array}
		\right.
		$$ 
		with $\Theta= \max \{\|\nabla \phi\|_{\infty}, \sqrt{\|V-E\|_{\infty}},1\}$. 
	\end{thm}
	
	We also adapt a gluing method develop by \cite*{Burq2022} to study Neumann and Dirichlet boundary condition, in order to study the localization on compact surface with semiclassical Robin boundary conditions: 
	
	\begin{thm}
		Assume that $M$ is a compact Riemannian surface with $W^{2,\infty}$ boundary. For $h$ small enough, let be  $u_h\in H^1(M,\R)$ weak solution of $$-\frac{h^2}{K_h}\mathrm{div}_g( K_h\nabla_g )u_h + V_h u_h =E_hu_h \mbox{ on } M,$$ where  $K_h=e^{\phi_h/h}$ with $\phi_h=\phi+\mathcal{O}_{W^{1,\infty}}(h)\in W^{1,\infty}(M,\R)$, $V_h=V+\mathcal{O}_{L^\infty}(h)\in L^\infty(M,\R)$ and $E_h= E +\mathcal{O}_{\R}(h) \in\R$. Moreover, assume $u_h$ verify the Robin boundary condition: $$h\,\alpha(x)\,\partial_n {u_h}_{|\partial M}(x)+\beta(x)\,{u_h}_{|\partial_M}(x)=0,$$ with $n$ the normal exterior to the boundary and $\alpha$ and $\beta$ two functions on the boundary such that 
		\begin{itemize}
			\item if $\alpha(x)=0$ for $x\in \partial M$, then $\alpha=0$ on the connected component of $\partial M$ that contained $x$,
			\vspace{1mm}
			\item $\frac{\beta}{\alpha}\in W^{1,\infty}(\partial M\setminus \{\alpha=0\})$,
			\vspace{1mm}
			\item $\{\alpha=0\}\cap\{\beta=0\}=\emptyset$.
		\end{itemize}
	
	Then for any  ball $B_r\in \overline{M}$ of radius $r<1/4$:
	
	\begin{equation}
	\sup\limits_{B_r}|u_h| \geq  r^{ L\left(\frac{\epsilon(h)+\mathrm{osc}(\phi_h)+\sup\left|\frac{\beta}{\alpha}\right|1_{\{\alpha \neq 0\}}}{h}\right)} \sup_{M}|u_h|
\end{equation}
where $L$ is an absolute constant, $\mathrm{osc}(\phi_h)= \max_M \phi_h -\min_M \phi_h$  and 
$$
\epsilon(h):=  \left\{
\begin{array}{ll}
	\Theta\log\left(\frac{\Theta}{h}\right) & \mbox{if } 	||\nabla \phi||_{\infty}+\| \frac{\beta}{\alpha} 1_{\{\alpha \neq 0\}}\|_{W^{1,\infty}}>0,\\
	\Theta\log\left(\frac{\Theta}{h}\right)^{1/2} & \mbox{else.}
\end{array}
\right.
$$ 
with $\Theta= \max \{\|\nabla \phi\|_{\infty}+\| \frac{\beta}{\alpha} 1_{\{\alpha \neq 0\}}\|_{W^{1,\infty}}, \sqrt{\|V-E\|_{\infty}},1\}$. 
\end{thm}

	The outline of this article is as follows. In Section~\ref{sect: Preliminaire} we first state and adapt to our framework some result from elliptic theory that we will use in the proof of Theorem \ref{thm: mainthm} (Section~\ref{sect: elleptic_theory}). We conclude the section~\ref{sect: Preliminaire} with a discussion on hole cuttings within the plane framework (Section~\ref{sect: decoupespace}).
	
	Section~\ref{sect: de SchrodAHarmonique} is dedicated to reduce the eigenfunctions to harmonic functions. This process unfolds in three steps: (1) reduction to a divergence-free equation (Section~\ref{subsect: eqdivsuresptroue}), (2) application of quasiconformal mappings theory to transform the divergence-free equation to an harmonic one (Section~\ref{subsect: chg de var quasi-conforme}), (3) analysis of the geometric impact of these mappings on our framework (Section~\ref{subsect: distortion de lespace}). To start this process, we will need to construct eigenfunctions of the weighted Schrödinger operator with some "good properties" (Section~\ref{sect: constrsolvarphi}).
	
	In Section~\ref{sect: end_proof_main}, we finish the proof of Theorem \ref{thm: mainthm} by making use of the Theorem $5.3$ on harmonic functions from \cite*{logunov2020landisconjectureexponentialdecay} result on harmonic.
	
	We finish the article by the study of eigenfunctions on compact surface: closed surface (Section~\ref{sect: closed_surface}), and compact surface with boundaries conditions (Section~\ref{sect: compact_surface_with_boundary}). Our way to deal with the Robin boundary condition will lead us to replace it with a Neumann condition at the cost of conjugating the Schrödinger operator with an exponential term.
	
	\textbf{Aknowledgement:} The author acknowledges the funding from the European Research Council (ERC) under the European Union’s Horizon 2020 research and innovation programme (Grant agreement 101097172 – GEOEDP). 
	
	\section{Preliminaries}
	\label{sect: Preliminaire}
	
	\subsection{Renormalization and notation}
	\label{sect: Renorm_and_not}
	
	In order to simplify some of the following calculus and result we will, renormalize the eigenvalue equation by $\Theta=\max\{\|V-E\|_{\infty},\|\nabla \phi\|_{\infty}^2,1\}$:
	  $$ \left(-\frac{h^2}{\Theta K_h}\mathrm{div}(K_h\nabla) +\frac{(V_h-E_h)}{\Theta}\right)u_h=0. $$
	So with some abuse of notation we will replace:
	\begin{itemize}
		\item $\frac{(V_h-E_h)}{\Theta}$ by $V_h= V +\mathcal{O}_{L^{\infty}}(h)$ with $\|V\|_{\infty}\leq 1$;
		\item $\frac{h}{\sqrt{\Theta}}$ by $h$;
		\item and so $\frac{\phi_h}{\sqrt{\Theta}}$ by $\phi_h= \phi + \mathcal{O}_{W^{1,\infty}}(h)$ with $\|\nabla \phi\|_{\infty}\leq 1$.
	\end{itemize}
	
	So we interest to the equation:
	\begin{equation}
		P_h u_h:= \left(-\frac{h^2}{K_h}\mathrm{div}(K_h\nabla) +V_h \right)u_h=0.
	\end{equation}
	Thus we will refer now the $E_h-$eigenfunction $u_h$ to a zero-eigenfunction of $P_h$.
	
	We assume $h$ small enough such that $\|V_h\|_{\infty}\leq 2$ and $\| \nabla \phi_h\|_{L^\infty}\leq 2$. 
	
	\subsection{Some elliptic theory results}
	\label{sect: elleptic_theory}
	
	In this section, we state some result from elliptic theory that we will use later. They are modifications of precedent lemmas from the article of \cite{logunov2020landisconjectureexponentialdecay} to correspond to our framework. Their proofs can be found in the Appendix section.
	
	\begin{lem}
		\label{lem: rég_C1}
		If $P_{h}u_{h}=0$ on an open set $\Omega$ of $\R^2$ in the sense of distributions and $u_{h}\in H_{loc}^1(\Omega)$, then $u_{h} \in \mathcal{C}^1(\Omega)$.
	\end{lem}
	This lemma shows that our weak solutions actually possess a $\mathcal{C}^1$ regularity and so are locally Lipschitz. This will be useful for reaching inequality and properties which are inaccessible with only $H^1$ estimate (see for instance the proof proposition \ref{prop: divtypeq}). 
	
	The next lemma concerns the existence of solutions for some Divergence type equations:
	\begin{lem}
		\label{lem: constsolSchrodsmallPconst}
		Let $\Omega_h$ be an open set of $B(0,R_0)$ with Poincaré constant $k(h)<<h$. 
		There is a weak solution (in the sense of distributions) $u\in H_0^1(\Omega_h,\R)$ to the equation
		\begin{equation}
			h^2K_h^{-1}\mathrm{div}(K_h\nabla u)=v
		\end{equation}
		for all $v\in L^2\cap L^\infty(\Omega_h,\R )$: for all $\chi \in \mathcal{C}_0^\infty(\Omega_h,\R)$ $$\int_{\Omega_h}h^2K_h\nabla u\cdot \nabla\chi=\int_{\Omega_h} K_hv\chi.$$
		The solution $u$ verify	
		\begin{itemize}
			\item $||u||_{H_{k(h)}^{1}(\Omega_h,\, K_h\mathrm{d}x)}\leq A\frac{k(h)^2}{h^2}||v||_{L^2(\Omega_h,\, K_h\mathrm{d}x)}$,
			\item $||u\, \sqrt{K_h}||_{L^\infty(\Omega_h)}\leq  A\frac{k(h)^2}{h^2}||v\, \sqrt{K_h}||_{L^\infty(\Omega_h)}$,
		\end{itemize}
		with $A>0$ an absolute constant and where $$||\cdot||_{H_{k(h)}^{1}(\Omega_h, K_h\mathrm{d}x)}:= \sqrt{\int_{\Omega_h} \left(|\cdot|^2+|k(h)\nabla \cdot|^2 \right)K_h\mathrm{d}x}.$$ 
	\end{lem}
	With this lemma, we will construct solutions $\varphi_h$ to equation $P_{h}\varphi_h=0$ with some adequate properties. Look at, for instance, the section \ref{sect: constrsolvarphi}.
	
	The following lemma will provide information about the structure of the nodal set of the Schrödinger eigenfunctions: 
	\begin{lem}
		\label{lem: defr0}
		Let $u_h$ be a $H^1-$solution of $P_{h}u_h=0$ on a ball $B(x,r_h)$ of $B(0,R_0)$ where $r_h<r_0h$ with $r_{0}$ a small constant independent of $h$ and $u_h$, and  $0<h\leq1$. If $u_h>0$ on $\partial B(x,r_h)$, then $u_h>0$ on $B(x,r_h)$. 
	\end{lem}
	We will use in fact the converse of this lemma to say that there is an universal radius $r_0$ for which at any point $x$ in the zero set of $u_h$, the circle with centre $x$ and radius $rh$ with $r\leq r_0$ intersects the zero-set of $u_h$ or the boundary of $B(0,R_0)$. This lemma will be use directly in the following section.

	\subsection{Perforation of space}
	\label{sect: decoupespace}
	
	Let’s finish the preliminaries by setting the stage for the proof. In this section, we modify our framework space by replacing $B(0,R_0)$ by a set with smaller Poincaré constant.  
	
	Let $F_{0,h}$ be the nodal set of $u_h$ inside $B(0,R_0)$.  By the lemma \ref{lem: defr0}, we dispose of a radius $r_0>0$ such that for $h$ small enough, for all $x\in F_{0,h}$ and for all $r<r_0$,  the circle $C(0,hr)$ intersect $F_{0,h}$ or the boundary of $B(0,R_0)$.
	Let $\epsilon (h) >0$ and $C>2$ be respectively a function of $h$ such that $\epsilon(h)\rightarrow 0$ when $h\rightarrow 0$ and a positive constant. These two parameters will be determined later. We assume that $h<<\epsilon(h)$ for $h$ small enough. 
	
	We begin by filling $B(0,R_0)\setminus  F_{0,h}$ with balls of radius $(C+1)h\epsilon(h)$.  In this configuration, all points of $B(0,R_0)$ are at a distance smaller than $2(C+1)h\epsilon(h)$ from a ball, from $ C(0,R_0)$ (the boundary of $B(0,R_0)$) or from $F_{0,h}$  (otherwise we could add one more ball at these points).
	
	We define $F_{1,h}$ as the set of balls with radius $h\epsilon(h)$, centre on the same points of the balls above. Thus all the balls of $F_{1,h}$ are separated by a distance bigger than $Ch\epsilon(h)$ from any other balls of $F_{1,h}$,  $C(0,R_0)$ and  $F_{0,h}$. 
	By construction, the set $F_{1,h}\cup F_{0,h}\cup C(0,R_0)$ is at least $3Ch\epsilon(h)-$dense in $B(0,R_0)$.
	
	From the set $F_{1,h}$ we remove the balls that are $Ch\epsilon(h)$ near from $0$ and $x_{\max,h}$ define by $|u_h(x_{\max,h})|=\sup_{B(0,R_0/2)}|u_h|$. Then the new set $\tilde{F}_{1,h}\cup F_{0,h}\cup C(0,R_0)$ is at least $10Ch\epsilon(h)-$dense in $B(0,R_0)$. To simplify the notation, we will continue to use $F_{1,h}$ as $\tilde{F}_{1,h}$.

	We suppose that $\epsilon(h)<r_0$ (which will be true for $h$ small enough).
	Thanks to this perforation, we get a nice control on the Poincaré constant of the set $$\Omega_h:=B(0,R_0)\setminus \textbf{(}  F_{1,h} \cup F_{0,h}\textbf{)}:$$
	
	\begin{lem}
		For $h$ small enough, the Poincaré constant $k(h)$ of $\Omega_h$, define here as the smallest constant such that $$\int |f|^2K_hdx\leq k(h)^2\int |\nabla f|^2K_hdx$$ for all $f \in \mathcal{C}_0^\infty(\Omega_h)$, is smaller than $C'h\epsilon(h)$ where $C'$ is constant only dependent of $C$.
	\end{lem}
	
	\begin{proof}
		
		Let $f \in \mathcal{C}_0^\infty(\Omega_h)$. We extend $f$ by $0$ on $B(0,R_0)$. Let $z \in F_{1,h}\cup F_{0,h}\cup C(0,R_0) $. For $h$ small enough, we know that for all $0<r_h<h\epsilon(h)<hr_0$, $$C(z,r_h)\cap \textbf{(}F_{1,h}\cup F_{0,h}\cup C(0,R_0)\textbf{)}\neq \emptyset.$$ So
		\begin{equation}
			\max_{C(z,r)} |f|\leq \int_{C(z,r)}|\nabla f|.
		\end{equation}  
		We use this to obtain this first Poincaré Inequality on $B(z,h\epsilon(h))$:
		\begin{align}
			\int_{B(z,h\epsilon(h))}& |f|^2 = \int_0^1 \left(\int_{C(z,r)}|f|^2\right) \mathrm{d}r\nn\\
			& \leq \int_0^1 |C(z,r)|\max_{C(z,r)} |f|^2 \mathrm{d}r \leq \int_0^1 |C(z,r)|\left(\int_{C(z,r)}|\nabla f|\right)^2 \mathrm{d}r\nn\\
			&\leq \int_0^1 |C(z,r)|^2\left(\int_{C(z,r)}|\nabla f|^2\right)\mathrm{d}r\leq C_1h^2\epsilon(h)^2 \int_{B(z,h\epsilon(h))}|\nabla f|^2.
		\end{align}
		
		Then we can find $r_h\in (\,h\epsilon(h)/2\, ;\, h\epsilon(h)\,)$ such that:
		
		\begin{equation}
			\label{eq: proofsmallPutil1}
			\int_{C(z,r_h)}|f|^2\leq \frac{C_2}{h\epsilon(h)} \int_{B(z,h\epsilon(h))} |f|^2\leq C_3h\epsilon(h) \int_{B(z,h\epsilon(h))}|\nabla f|^2.
		\end{equation}
		
		Let be $\Gamma_\theta$ the segment starting at the point $x_\theta=z+r_he^{i\theta}$ and ending at the point $z+10Ch\epsilon(h)e^{i\theta}$.
		We have:
		\begin{multline}
			\max_{x\in \Gamma_\theta} |f|^2\leq\left(|f(x_\theta)|+\int_{\Gamma_\theta} |\nabla f| \right)^2\leq 2 \left(|f(x_\theta)|^2+\left(\int_{\Gamma_\theta} |\nabla f|\right)^2 \right)\nn\\
			\leq  2 \left(|f(x_\theta)|^2+|\Gamma_\theta|\int_{\Gamma_\theta} |\nabla f|^2 \right)\leq 2 \left(|f(x_\theta)|^2+C_4 h\epsilon(h)\int_{\Gamma_\theta} |\nabla f|^2 \right).
		\end{multline}
		
		Therefore we can extend the last Poincaré inequality to the ball $B(z,3Ch\epsilon(h))$:
		\begin{align}
			\int_{B(z,3Ch\epsilon(h))} &|f|^2= \int_{B(z,h\epsilon(h))} |f|^2+\int_{B(z,3Ch\epsilon(h))\setminus B(z,h\epsilon(h))} |f|^2 \nn\\
			&\leq C_1h^2\epsilon(h)^2 \int_{B(z,h\epsilon(h))}|\nabla f|^2 +(10Ch\epsilon(h))^2\int_0^{2\pi}\max_{x\in \Gamma_\theta} |f|^2\nn\\
			&\leq C_1h^2\epsilon(h)^2 \int_{B(z,h\epsilon(h))}|\nabla f|^2 +2(10Ch\epsilon(h))^2 \int_0^{2\pi}\left(|f(x_\theta)|^2+C_4 h\epsilon(h)\int_{\Gamma_\theta} |\nabla f|^2 \right).
		\end{align}
		And thank to the equation \ref{eq: proofsmallPutil1}, we conclude:
		\begin{equation}
				\int_{B(z,3Ch\epsilon(h))} |f|^2\leq C_5 h^2\epsilon(h)^2 \int_{B(z,3Ch\epsilon(h))} |\nabla f|^2.
		\end{equation} 
		
		Because $B(0,R_0)$ is bounded and $ F_{1,h}\cup F_{0,h}\cup C(0,R_0)$ is $10Ch\epsilon(h)$ dense in it, we can choose a finite collection of point  such that $B(z,10Ch\epsilon(h))$ cover $B(0,R_0)$ and each point is covered a bounded number of times.
		
			Consequently, we have that for any $f\in \mathcal{C}_0^\infty(\Omega_h)$
		$$||f||_{L^2(\Omega_h)}\leq \tilde{C}h\epsilon(h)||\nabla f||_{L^2(\Omega_h)}.$$ If we replace by density the function $f$ by $f\,K_h=f\,e^{\phi_h/h}\in H_0^1(\Omega_h)$, we can see that $$||f\,K_h||_{L^2(\Omega_h)}\leq \tilde{C}h\epsilon(h)||\nabla f\,K_h||_{L^2(\Omega_h)}\leq \tilde{C}h\epsilon(h)||(\nabla f)K_h||_{L^2(\Omega_h)}+\tilde{C}\epsilon(h)||\nabla \phi_h||_\infty||f\,K_h||_{L^2(\Omega_h)}$$ with $||\nabla \phi_h||_\infty$ bounded independently of $h$. Thus, for $\epsilon(h)$ small enough, we get the Poincaré estimate.
	\end{proof}

	In the following, the Poincaré constant will always be small. We can expect $C$ to be very large, so $C'$ will be too. this imply that $h\epsilon(h)$ is small.
	
	\section{From Schrödinger to Harmonic equation}
	\label{sect: de SchrodAHarmonique}
	
	The aim of this section is to substitute our study of zero-eigenfunctions of the operator $P_{h}$ by a study of harmonic functions which are simpler to manipulate. For that, we'll first replace our functions by free divergence ones. Then we'll use the theory of quasi-conformal mappings to reduce to the case of harmonic functions. 
	
	\subsection{Construction of solutions to the zero-Eigenvalue equation of $P_{h}$}
	\label{sect: constrsolvarphi}
	In order to reduce to a divergence-type equation, we will modify the eigenfunction $u_h$ thanks to an another function that we will construct in this section. For that we will use the smallness of $\Omega_h$'s Poincaré constant (see the precedent section \ref{sect: decoupespace}). We will first show the following proposition:
	
	\begin{prop}
		\label{prop: constrsmallsolSchrod}
		Let $\Omega\subset B(0,R_0)$ an open set with Poincaré constant $k(h)$.
		If $k(h)<<h$, then there is a non-zero solution $\varphi_{h}\in H_{loc}^1(\Omega)$ to the equation $P_{h}\varphi_h=0$ on $\Omega$ in the sense of distributions such that  $\varphi_{h}=\tilde{K}_h^{-1/2}+\tilde{\varphi}_h$ with:
		
		\begin{itemize}
			\item $\tilde{K}_h\in \mathcal{C}^\infty(\overline{B(0,R_0)},\R_{+}^{*})$, such that 
			\begin{equation}
				\label{eq: controlsmallsol_eq_1}
				|K_h\tilde{K}_h^{-1}-1|\lesssim \frac{k(h)}{h}||\nabla \phi_h||_{\infty}+\mathcal{O}\left(\frac{k(h)^2}{h^2}\right);
			\end{equation}
			\item $\tilde{\varphi}_h \in H_0^{1}(\Omega, \R)$, such that 
			\begin{align}
				\label{eq: controlsmallsol_eq_2}
				&||\tilde{\varphi}_h||_{H_{k(h)}^{1}(\Omega,\,K_h\mathrm{d}x)} \lesssim \frac{k(h)}{h}||\nabla \phi_h||_{\infty}+\mathcal{O}\left(\frac{k(h)^2}{h^2}\right),\\
				\label{eq: controlsmallsol_eq_3}
				&||\tilde{\varphi}_h\sqrt{K_h}||_{\infty}\lesssim \frac{k(h)}{h}||\nabla \phi_h||_{\infty}+\mathcal{O}\left(\frac{k(h)^2}{h^2}\right);
			\end{align}
		\end{itemize}   
		where $$||\cdot||_{H_{k(h)}^{1}(K_h\mathrm{d}x)}:= \sqrt{\int \left(|\cdot|^2+|k(h)\nabla \cdot|^2 \right)K_h\mathrm{d}x}.$$
		Consequently, for $h$ small enough, 
		\begin{equation}
			||1-\sqrt{K_h}\varphi_h||_\infty \lesssim \frac{k(h)}{h} ||\nabla \phi_h||_{\infty}+\mathcal{O}\left(\frac{k(h)^2}{h^2}\right). 
		\end{equation}
	\end{prop} 
	
	\begin{proof}
		Our goal is to apply recursively the lemma \ref{lem: constsolSchrodsmallPconst} to construct the function  $\tilde{\varphi}_{h}$. For that will need to construct the function $\tilde{K}_h$ first such that  $||K_h^{1/2}P_{h}\tilde{K}_h^{-1/2}||_\infty$ don't diverge in $h$ too fast.
	Let be $\rho \in \mathcal{C}_0^{\infty}(B(0,1))$ such that $\int \rho=1$ and $0\leq \rho\leq 1$. We use this function as a smoothing kernel to define the regularisation of $\phi_h$ on $B(0,R_0)$: 
	\begin{equation}
		\tilde{\phi}_h(x):= \int \frac{1}{k(h)^2}\rho \left(\frac{t}{k(h)}\right)\phi_h(x-t)dt, \quad \mbox{ for } x\in B(0,R_0).
	\end{equation}
	which is well defined for $k(h)$ small enough small enough because $\phi_h\in W^{1,\infty}(B(0,R))$. Moreover we have $$|\tilde{\phi}_h-\phi_h|(x)\leq k(h) \|\nabla \phi_h\|_\infty,\,|\nabla \tilde{\phi}|(x)\leq \|\nabla \phi_h\|_\infty\,\mbox{ and } |\Delta \tilde{\phi}_h |(x)\lesssim k(h)^{-1}\|\nabla \phi_h\|_\infty.$$  Let note $\tilde{K}_h=e^{\tilde{\phi}_h/h}$.
	We remind that with the supposition on $h$ from the section \ref{sect: Renorm_and_not}, $\|\nabla \phi_h\|_\infty\leq 2$ and $\|V_h\|_\infty\leq 2$.
	
	Firstly we have 
	\begin{multline}
		\left|1-K_h\tilde{K}_h^{-1} \right|\leq \left|1-e^{(\phi_h-\tilde{\phi}_h)/h}\right|\leq 
		\left|1-e^{- \|\nabla \phi_h\|_\infty\,k(h)/h}\right|\lesssim \frac{k(h)}{h}\|\nabla \phi_h\|_\infty+\mathcal{O}\left(\frac{k(h)^2}{h^2}\right).
	\end{multline}
	Which imply the convergence of $K_h/\tilde{K}_h$ to $1$ when $h\to 0$ and the estimate \ref{eq: controlsmallsol_eq_1}.
	
	Secondly, we check the $h$ dependency of $||K_h^{1/2}P_{h}\tilde{K}_h^{-1/2}||_\infty$:
	\begin{multline}
		\left|K_h^{1/2}P_{h}\tilde{K}_h^{-1/2}\right| = \left|-h^2K_h^{1/2}\Delta \tilde{K}_h^{-1/2}- h^2K_h^{-1/2}\nabla K_h \nabla \tilde{K}_h^{-1/2} + V_h \left(K_h\tilde{K}_h^{-1}\right)^{1/2}  \right|.
	\end{multline}
	
	Firstly, with what we see above, thus
	\begin{equation}
	 |V_h| \left(K_h\tilde{K}_h^{-1}\right)^{1/2}\lesssim ||V_h||_\infty \left(1+ \frac{k(h)}{h}\|\nabla \phi_h\|_\infty+\mathcal{O}\left(\frac{k(h)^2}{h^2}\right)\right)
	\end{equation}
 for   $k(h)/h$ small enough. Secondly
	$$-h^2K_h^{-1/2}\nabla K_h \nabla \tilde{K}_h^{-1/2} = \frac{1}{2}\left(K_h\tilde{K}_h^{-1}\right)^{1/2} \nabla \phi_h \nabla \tilde{\phi}_h$$
	and by construction $||\nabla \tilde{\phi}_h||_\infty\leq ||\nabla \phi_h||_\infty$.  So
	\begin{equation}
		\left|h^2K_h^{-1/2}\nabla K_h \nabla \tilde{K}_h^{-1/2}\right|  \lesssim ||\nabla \phi_h||_\infty^2 \left(1+ \frac{k(h)}{h}\|\nabla \phi_h\|_\infty+\mathcal{O}\left(\frac{k(h)^2}{h^2}\right)\right)
	\end{equation}
	 for all $k(h)/h$  small enough. Finally
$$-h^2K_h^{1/2}\Delta \tilde{K}_h^{-1/2}=\frac{1}{2}\left(K_h\tilde{K}_h^{-1}\right)^{1/2}(h\Delta \tilde{\phi}_h-|\nabla \tilde{\phi}_h|^2)$$
	and $|| h\Delta \tilde{\phi}_h||_\infty\lesssim ||\nabla \phi_h||_\infty\, h/k(h)$. So
		\begin{multline}
		\left|h^2K_h^{1/2}\Delta \tilde{K}_h^{-1/2}\right| \lesssim 
		 \left(\frac{h}{k(h)} ||\nabla \phi_h||_\infty+||\nabla \phi_h||_\infty^2\right) \left(1+ \frac{k(h)}{h}\|\nabla \phi_h\|_\infty+\mathcal{O}\left(\frac{k(h)^2}{h^2}\right)\right)
		\end{multline}
	for all  $k(h)/h$  small enough.
	Hence for $k(h)/h$ small enough:
	\begin{multline}
		||K_h^{1/2}P_{h}\tilde{K}_h^{-1/2}||_{L^\infty(B(0,R_0))} \lesssim \\ \left(\frac{h}{k(h)}||\nabla \phi_h||_\infty +||\nabla \phi_h||_\infty^2+||V_h||_\infty \right)\left(1+ \frac{k(h)}{h}\|\nabla \phi_h\|_\infty+\mathcal{O}\left(\frac{k(h)^2}{h^2}\right)\right) \\
		\lesssim \frac{h}{k(h)}||\nabla \phi_h||_\infty  + \mathcal{O}\left(1\right).
	\end{multline}
	
	Now, with what we just check, we can apply the lemma \ref{lem: constsolSchrodsmallPconst} for $v=P_{h}\tilde{K}_h^{-1/2}\in L^2 \cap L^\infty(\Omega)$. We dispose then of a function $\varphi_{1,h}\in H_0^{1}(\Omega)$ weak solution of $$\frac{h^2}{K_h}\mathrm{div}(K_h\nabla  \varphi_{1,h})=P_{h}\tilde{K}_h^{-1/2}$$ on $\Omega_h$ with 
	$$||\varphi_{1,h}||_{H_{k(h)}^{1}(\Omega,\, K_h\mathrm{d}x)} \leq A\frac{k(h)^2}{h^2}	||K_h^{1/2}P_{h}\tilde{K}_h^{-1/2}||_{L^2(B(0,R_0))}\lesssim  \frac{k(h)}{h}||\nabla \phi_h||_\infty\,+ \mathcal{O}\left(\frac{k(h)^2}{h^2}\right)$$
	and
	$$ ||\varphi_{1,h}\sqrt{K_h}||_{L^\infty(\Omega)}\leq  A\frac{k(h)^2}{h^2}	||K_h^{1/2}P_{h}\tilde{K}_h^{-1/2}||_{L^\infty(B(0,R_0))}\lesssim \frac{k(h)}{h}||\nabla \phi_h||_\infty\,+ \mathcal{O}\left(\frac{k(h)^2}{h^2}\right)$$
	with $A>0$ absolute constant from the lemma \ref{lem: constsolSchrodsmallPconst}.
	We can then repeat the processes with $\varphi_{2,h}\in H_0^1(\Omega)$ the function give by lemma \ref{lem: constsolSchrodsmallPconst} for $v=V_h\varphi_{1,h}$. We obtain then
	$$||\varphi_{2,h}||_{H_{k(h)}^{1}(\Omega,\, K_h\mathrm{d}x)}\leq A\frac{k(h)^2}{h^2}||V_h||_{\infty}||\varphi_{1,h}||_{H_{k(h)}^{1}(\Omega_h,\, K_h\mathrm{d}x)} $$
	and
	$$ ||\varphi_{2,h}\sqrt{K_h}||_{L^\infty(\Omega_h)}\leq  A\frac{k(h)^2}{h^2}||V_h||_\infty||\varphi_{1,h}\sqrt{K_h}||_{L^\infty(\Omega_h)}.$$
	
	Hence, by recurrence for $j\in \N^{*}$:
	$$||\varphi_{j,h}||_{H_{k(h)}^{1}(\Omega,\, K_h\mathrm{d}x)}  \lesssim \left(A||V_h||_\infty\, \frac{k(h)^2}{h^2} \right)^{j-1} \left(\, \frac{k(h)}{h} ||\nabla \phi_h||_\infty\, + \mathcal{O}\left(\frac{k(h)^2}{h^2}\right) \,\right),$$
	and
	$$ ||\varphi_{j,h}\sqrt{K_h}||_{L^\infty(\Omega_h)}\lesssim \left(A||V_h||_\infty\, \frac{k(h)^2}{h^2}\right)^{j-1}\left(\, \frac{k(h)}{h} ||\nabla \phi_h||_\infty\, + \mathcal{O}\left(\frac{k(h)^2}{h^2}\right) \,\right).$$
	
	So for $$A||V_h||_\infty\frac{k(h)^2}{h^2}\leq \frac{1}{2}$$ (which is possible for $k(h)^2/h^2\leq (2A)^{-1}$) the sum $\sum_j \varphi_{j,h}$ converge in $H_{k(h), 0}^{1}\cap L^\infty$ to $\tilde{\varphi}_{h}$. Moreover for $\chi \in \mathcal{C}_0^\infty(\Omega_h)$:
	
	\begin{multline}
		\int_{\Omega_h} h^2K_h \nabla \varphi_h\cdot\nabla \chi +K_hV_h\varphi_h\chi=\\
		 \sum_{j\geq 1} \int_{\Omega_h} h^2K_h \nabla \varphi_{j+1,h}\cdot\nabla \chi +K_hV_h\varphi_{j,h}\chi + \int_{\Omega_h}  h^2K_h \nabla \varphi_{1,h}\cdot\nabla \chi +\chi P_{h}\tilde{K}_h^{-1/2}=0.
	\end{multline}
	Thus  $P_{h}\varphi_h=0$ on $\Omega$.

\end{proof}

We apply this proposition to $\Omega= \Omega_h$. We dispose then of $\varphi_{h}\in H_{loc}^1(\Omega_h)$, solution of the equation $P_{h}\varphi_h=0$, such that 
$\varphi_h=\tilde{K}_h^{-1/2}+\tilde{\varphi}_h$ with $\tilde{K}_h^{-1/2}\in \mathcal{C}^\infty(\Omega_h)$ and $\tilde{\varphi}_h\in H_0^1(\Omega_h)$. Moreover, with the choices of renormalization from the section \ref{sect: Renorm_and_not} and because we assume $h<< \epsilon(h)$, for $h$ small enough
\begin{equation}
	\label{eq: convsmallsol1}
	||1-K_h^{1/2}\tilde{K}_h^{-1/2}||_\infty \lesssim \epsilon(h) \delta_{\|\nabla\|_{\infty}\neq 0}+ \mathcal{O}(\epsilon(h)^2),
\end{equation}
\begin{equation}
	\label{eq: convsmallsol2}
	||K_h^{1/2}\tilde{\varphi}_h||_\infty \lesssim \epsilon(h) \delta_{\|\nabla\|_{\infty}\neq 0}+ \mathcal{O}(\epsilon(h)^2)
\end{equation}
and
\begin{equation}
		\label{eq: convsmallsol3}
	||1-\sqrt{K_h}\varphi_h||_\infty \lesssim \epsilon(h) \delta_{\|\nabla\|_{\infty}\neq 0}+\mathcal{O}(\epsilon(h)^2). 
\end{equation}

	We extend $\tilde{K}_h^{-1/2}$ by $\tilde{K}_h^{-1/2}$ on $B(0,R_0)$ and we extend $\tilde{\varphi}_h$ by $0$ on $B(0,R_0)$.
	By the equations \ref{eq: convsmallsol1} and \ref{eq: convsmallsol2}, we can see that for $\epsilon(h)$ small enough $$1-\frac{1}{3}\leq K_h^{1/2}\tilde{K}_h^{-1/2}\leq1+ \frac{1}{3}\mbox{ and } ||\sqrt{K_h}\tilde{\varphi}_h||_\infty\leq \frac{1}{3}.$$ 
	 
	\begin{prop}
		\label{prop:deffh}
		$\varphi_h, \frac{1}{\varphi_h} \mbox{ and }\frac{u_h}{\varphi_h}\in H_{loc}^{1}(B(0,R_0))$ with
		 $$\nabla \frac{1}{\varphi_h}= -\frac{\nabla \varphi_h}{\varphi_h^2}, \nabla \varphi_h= \nabla {K}_h^{-1/2}+\nabla \tilde{\varphi}_h 1_{\Omega_h} $$
		 and
		 \begin{equation}
		 	\label{eq: toolusurphi}
		 	\nabla \frac{u_h}{\varphi_h}=\frac{\nabla u_h}{\varphi_h}-u_h\frac{\nabla \varphi_h}{\varphi_h^2}.
		 \end{equation} 
	\end{prop}
	
	\begin{proof}
		This proof is essentially the same that we can find in the facts $4.2$ and $6.12$ of \cite{logunov2020landisconjectureexponentialdecay}.
		
		\textbf{Formal check of $\varphi_h$ and $\varphi_h^{-1}$:}
		
		Firstly, let be $\psi_k\in \mathcal{C}_0^\infty(\Omega_h)$ a family of smooth test function such that $\psi_k \to \tilde{\varphi}_h$ in $H^1(\Omega_h)=H^1(\Omega_h,K_h\mathrm{d}x)$ (because $K_h$ and $K_h^{-1}\in L^\infty$) and  $||\psi_k\sqrt{K_h}||_\infty \leq 1/2$.
		We can extend $\psi_k$ to $B(0,R_0)$ by $0$ as a smooth function. Hence $\nabla \psi_k=\nabla \psi_k 1_{\Omega_h}$. Consequently $\psi_k \to \tilde{\varphi}_h 1_{\Omega_h}$ in $L^2(B(0,R_0))$ and a.e, and $\nabla\psi_k \to \nabla(\tilde{\varphi}_h) 1_{\Omega_h}$ in $L^2(B(0,R_0))$. So directly we have that $\varphi_h\in H_{loc}^1(B(0,R_0))$ with
		
		\begin{equation}
		\nabla \varphi_h= \nabla \tilde{K}_h^{-1/2}+\nabla \tilde{\varphi}_h 1_{\Omega_h}.
		\end{equation}
		
		Next, note that in $\Omega_h$:
		\begin{equation}
			\label{eq: proofH1loc1}
			\left|\frac{1}{\sqrt{K_h}(\tilde{K}_h^{-1/2}+\psi_k)}-\frac{1}{\sqrt{K_h}(\tilde{K}_h^{-1/2}+\tilde{\varphi}_{h})} \right|=\left|\frac{(\tilde{\varphi}_{h}-\psi_k)\sqrt{K_h}}{K_h(\tilde{K}_h^{-1/2}+\psi_k)(\tilde{K}_h^{-1/2}+\tilde{\varphi}_h)} \right|\lesssim |(\tilde{\varphi}_{h}-\psi_{k})\sqrt{K_h}|
		\end{equation}
		and
		\begin{align}
			\label{eq: proofH1loc2}
		&\left|\frac{\nabla (\tilde{K_h}^{-1/2}+\psi_k)}{\sqrt{K_h}(\tilde{K}_h^{-1/2}+\psi_k)^2}-\frac{\nabla(\tilde{K}_h^{-1/2}+\tilde{\varphi}_h)}{\sqrt{K_h}(\tilde{K}_h^{-1/2}+\tilde{\varphi}_h)^2} \right|\leq 
		|(\nabla \psi_k-\nabla\tilde{\varphi}_h)\sqrt{K_h}|\frac{1}{K_h(\tilde{K}_h^{-1/2}+\psi_k)^2} \nn\\
		&+|\nabla\tilde{\varphi}_h\sqrt{K_h}|\left|\frac{1}{K_h(\tilde{K}_h^{-1/2}+\psi_k)^2}-\frac{1}{K_h(\tilde{K}_h^{-1/2}+\varphi_h)^2}\right|.
		\end{align}
		
		It is clear that, by the extension of $\psi_k$, $$(\tilde{K}_h^{-1/2}+\psi_k)^{-1}=\tilde{K}_h^{1/2} \mbox{ and }\nabla (\tilde{K}_h^{-1/2}+\psi_k)^{-1}=\nabla \tilde{K}_h^{1/2} $$ on $B(0,R_0)\setminus \Omega_h$. Hence by the above observation \ref{eq: proofH1loc1},
		\begin{equation}
		\frac{1}{\sqrt{K_h}(\tilde{K}_h^{-1/2}+\psi_k)}\to \frac{1}{\sqrt{K_h}(\tilde{K}_h^{-1/2}+\tilde{\varphi}_h)}1_{\Omega_h}+ \sqrt{\frac{\tilde{K}_h}{K_h}}1_{B(0,R_0)\setminus \Omega_h}
		\end{equation}
		 in $L^2(B(0,R_0))$ and a.e. 
		 
		 We would like to show that
		 \begin{equation}
		 	\frac{1}{\sqrt{K_h}}\nabla\frac{1}{\tilde{K}_h^{-1/2}+\psi_k}\to -\frac{\sqrt{K_h}\nabla (\tilde{K}_h^{-1/2}+\tilde{\varphi}_h)}{K_h(\tilde{K}_h^{-1/2}+\tilde{\varphi}_h)^2}1_{\Omega_h}+ \frac{\nabla\sqrt{\tilde{K}_h}}{\sqrt{K_h}}1_{B(0,R_0)\setminus \Omega_h}.
		 \end{equation} 
		 
		By the above observation \ref{eq: proofH1loc2}, we have
		$$\int_{\Omega_h}|(\nabla \psi_k-\nabla\tilde{\varphi}_h)\sqrt{K_h}|^2\frac{1}{K_h^2(\tilde{K}_h^{-1/2}+\psi_k)^4}\lesssim \int_{\Omega_h} |(\nabla \psi_k-\nabla\tilde{\varphi}_h)\sqrt{K_h}|^2\to 0,$$
		
		and the function $$|\nabla\tilde{\varphi}_h\sqrt{K_h}|^2\left|\frac{1}{K_h(\tilde{K}_h^{-1/2}+\psi_k)^2}-\frac{1}{K_h(\tilde{K}_h^{-1/2}+\varphi_h)^2}\right|^2$$
		
		 have an integrable majorant $C|\nabla\tilde{\varphi}_h\sqrt{K_h}|^2$ and so tend to $0$ in $\Omega_h$.
		Thus we have that $(\tilde{K}_h^{-1/2}+\psi_k)^{-1} \to (\tilde{K}_h^{-1/2}+\tilde{\varphi}_h)^{-1}1_{\Omega_h}+ \tilde{K}_h^{1/2}1_{B(0,R_0)\setminus \Omega_h}$ in $H^1(B(0,R_0),K_h^{-1}\mathrm{d}x)=H^1(B(0,R_0))$. Thus $\varphi_h^{-1} \in H_{loc}^1(B(0,R_0))$ with 
		
		\begin{equation}
			\nabla \frac{1}{\varphi_h}=-\frac{\nabla (\tilde{K}_h^{-1/2}+\tilde{\varphi}_h)}{(\tilde{K}_h^{-1/2}+\tilde{\varphi}_h)^2}1_{\Omega_h}+ \nabla\tilde{K}_h^{1/2}1_{B(0,R_0)\setminus \Omega_h}=-\frac{\nabla\tilde{K}_h^{-1/2}+\nabla\tilde{\varphi}_h1_{\Omega_h}}{(\tilde{K}_h^{-1/2}+\tilde{\varphi}_h1_{\Omega_h})^2}= -\frac{\nabla \varphi_h}{\varphi_h^2}.
		\end{equation}

		\textbf{formal check of $\frac{u_h}{\varphi_h}$:} 
		
		Finally fo the last function $\frac{u_h}{\varphi_h}$, the result come from the fact that $u_h$ and $(\varphi_{h})^{-1}$ are in $H_{loc}^{1}\cap L_{loc}^\infty$. Thus we can define their product in $H_{loc}^1(B(0,R_0))$ with
		\begin{equation}
			\nabla \frac{u_h}{\varphi_h}=\frac{\nabla u_h}{\varphi_h}-u_h\frac{\nabla \varphi_h}{\varphi_h^2}.
		\end{equation} 
		
		The fact is purely local, so we can work on $ B(0,r)\subset B(0,R_0)$ and $u_h, \varphi_h^{-1} \in H^{1}(B(0,r))\cap L^\infty(B(0,r))$.
		
		Let be $\delta_0>0$ small and
		let $$\psi_\delta(t)=\frac{1}{\delta^2}\psi(t/\delta)$$ be a smoothing kernel with $\delta<\delta_0$, $Supp (\psi_\delta)\subset B(0,\delta)$ and $\int \psi_{\delta}=1$.
		Then $u_{\delta}=\psi_\delta* u_h$ and $v_\delta=\psi_\delta*\varphi_h^{-1}$
		converge to $u_h$ and $\varphi_h^{-1}$ in $L^2(B(0,r-\delta_0))$ and a.e. in $B(x,r/2)$ as $\delta \to 0$.  Moreover we have $|u_\delta|< \|u_h\|_{L^\infty(B(0,r))}$, $|v_{\delta}|< \|\varphi_h^{-1}\|_{L^\infty(B(0,r))}$ in $B(0,r-\delta_0)$.
	
		We know that the convergence of $u_\delta v_\delta$ holds  a.e. in $B(0,r-\delta_0)$ and
		$|u_\delta v_\delta|$ are bounded by $\|u_h\|_\infty\|\varphi_h^{-1}\|_\infty$ a.e. in $B(0,r-\delta_0)$. So by the Lebesgue dominated convergence Theorem $u_\delta v_\delta \to u_h\varphi_h^{-1}$  in  $L^2(B(0,r-\delta_0))$.
		
		We want to show that the equation \ref{eq: toolusurphi} hold in the sense of $H^{1}(B(0,r-\delta_0))$.
		
		 It is clear that the right term is in $L^2(B(0,r-\delta_0))$ because $u_h$ and $\varphi_h^{-1}$ are bounded and their gradients are in $L^2(B(0,r-\delta_0))$.
		
		Let be $\eta \in C_0^\infty(B(0,r-\delta_0))$ a test function. We have
		$$ \int \frac{u_h}{\varphi_h} \nabla \eta = \lim_{\delta \to 0} \int u_\delta v_\delta \nabla \eta = - \lim_{\delta \to 0}  \int (\nabla u_\delta v_\delta + u_\delta \nabla v_\varepsilon ) \eta$$
		$$ =- \lim_{\delta \to 0} \left[ \int \left((\nabla u_\delta - \nabla u_h) v_\delta + u_\delta (\nabla v_\delta - \nabla \varphi_h^{-1}) \right) \eta + \int ( \nabla u_h v_\delta + u_\delta \nabla  \varphi_h^{-1} ) \eta  \right].$$
		
		Note that $$\int (\nabla u_\delta - \nabla u_h) v_\delta \eta  \to 0$$ because $\nabla u_\delta\to  \nabla u_h$ in $L^2(B(0,r-\delta_0))$ and  $|v_\delta \eta | < \|\varphi_h^{-1}\|_\infty\| \eta\|_\infty$ in  $B(0,r-\delta_0)$. Similarly, $$ \int u_\delta (\nabla v_\delta - \nabla \varphi_h^{-1}) \eta \to 0.$$
		Finally, by the Lebesgue dominated convergence Theorem $$ \int ( \nabla u_h v_\delta + u_\delta \nabla  \varphi_h^{-1} ) \eta  \to \int (\varphi_h^{-1} \nabla u_h + u_h \nabla  \varphi_h^{-1} ) \eta$$
		because the convergence of the functions holds a.e. in  $B(0,r-\delta_0)$ and there is the integrable majorant
		$(|\nabla u_h| \|\varphi_h^{-1}\|_\infty + \|u_h\|_\infty |\nabla  \varphi_h^{-1} |) \|\eta\|_\infty$.
		Thus 
		$$\int\frac{u_h}{\varphi_h} \nabla \eta = - \int ( \nabla u_h \varphi_h^{-1} + u_h \nabla  \varphi_h^{-1} ) \eta.$$
		Hence $\frac{u_h}{\varphi_h}\in H^1(B(0,r-\delta_0))$ for all $r<R_0$ and $\delta_0$. So $\frac{u_h}{\varphi_h}\in H_{loc}^1(B(0,R_0))$.
	\end{proof}
	
		We can now define the new function $f_h$ that will replace $u_h$ as
		\begin{equation}
			f_h:=\frac{u_h}{\varphi_h}\in H_{loc}^1(B(0,R_0)).
		\end{equation} 
	
	\subsection{Divergence type equations on the perforated space}
	\label{subsect: eqdivsuresptroue}
	
	In this section we show that the function $f_h=\frac{u_h}{\varphi_h}\in H_{loc}^1(B(0,R_0))$ is solution of a Divergence type equation on the space $\Omega'_h:=B(0,R_0)\setminus F_{1,h} $.
	
	\begin{prop}
		\label{prop: divtypeq}
		The function $f_h$ is a weak solution  of
		\begin{equation}
			\label{eq: divtypeq}
			-h^2\mathrm{div}(K_h\varphi_h^2\nabla f_h)=0
		\end{equation}
		on $\Omega'_h$.
	\end{prop}
	
	\begin{proof}
		
		This is due to the facts that $u_h$ and $\varphi_h$ are weak solutions of the same equation and that $u_h=0$ on $F_{0,h}$.
		
		Indeed, let $\chi_h \in C_0^\infty(\Omega'_h)$, the divergence equation \ref{eq: divtypeq} rewrite as
		$$ \int_{\Omega'_h} h^2 K_h\varphi_h^2 \nabla \left(\frac{u_h}{\varphi_h} \right) \nabla \chi_h =0.$$ By proposition \ref{prop:deffh}, we know that $$\nabla \frac{u_h}{\varphi_h}=\frac{\nabla u_h}{\varphi_h}-u_h\frac{\nabla \varphi_h}{\varphi_h^2}$$  in $B(0,R_0)$.
		So we need to prove
		
		\begin{equation} 
			\label{eq: toproove}
			\int_{\Omega'_h} h^2K_h(\varphi_h \nabla u_h - u_h \nabla \varphi_h)  \cdot  \nabla \chi_h =0.
		\end{equation}
		
		Since $u_h$ is a solution to $P_{h}u_h =0$ in $B(0,R_0)$, we have 
		\begin{equation} 
			\label{eq: divutil1}
			\int_{\Omega'_h} h^2K_h\nabla u_h \cdot (\varphi_h \nabla \chi_h + \chi_h \nabla \varphi_h) = \int_{\Omega'_h}K_h V_h\varphi_h u_h \chi_h
		\end{equation}
		(this equality is usually true for $\varphi_h$ smooth. We obtain it by using density of smooth functions in $ H_{loc}^1(B(0,R_0))$).
		
		Consider a function $\eta_\delta \in C_0^\infty(B(0,R_0) \setminus F_{0,h})$ such that $0\leq \eta_\delta\leq1$, $\eta_h=1$ outside a $\delta-$neighbourhood of $F_{0,h}\cup C(0,R_0)$ and descends from $1$ to $0$ in the $\delta-$neighbourhood with $|\nabla \eta_\delta| < A_h/\delta$.
		
		Since $P_{h}\varphi_h =0$ in $\Omega_h$ in the sense of distribution and $u_h\chi_h\eta_\delta \in C_0^1(\Omega_h)$ (because $u_h \in \mathcal{C}^1$ by lemma \ref{lem: rég_C1}), we have
		
		\begin{equation}
			\label{eq: divutil2}
			\int_{\Omega_h} h^2K_h \nabla \varphi_h \cdot(\chi_h \nabla u_h \eta_\delta +u_h\nabla \chi_h \eta_\delta + u_h \chi_h \nabla \eta_\delta)=  \int_{\Omega_h} K_hV_h\varphi_h u_h \chi_h\eta_\delta.
		\end{equation}
		
		Note that $$\int_{\Omega_h} K_hV_h\varphi_h u_h \chi_h\eta_\delta \rightarrow \int_{\Omega_h} K_hV_h\varphi_h u_h \chi_h  \mbox{ as } \delta \to 0$$ (by Lebesgue dominated convergence in $\Omega_h$ and because $\eta_\delta \to 1$ a.e in $\Omega_h$). 
		Note also that
		$\chi_h (\nabla u_h) \eta_\delta \to \chi_h\nabla u_h \mbox{ a.e in } \Omega_h \mbox{ and } u_h(\nabla \chi_h) \eta_\delta \to  u_h\nabla \chi_h \mbox{ a.e in } \Omega_h. $
		Hence
		\begin{equation} 
			\label{eq: divutil3}
			\int_{\Omega_h}h^2K_h \nabla \varphi_h \cdot(\chi_h \nabla u_h \eta_\delta +u_h\nabla \chi_h \psi_h) \to \int_{\Omega_h} h^2K_h  \nabla \varphi_h \cdot (\chi_h\nabla u_h + u_h\nabla \chi_h)  \mbox { as } \delta \to 0
		\end{equation}
		by the Lebesgue dominated convergence Theorem.
		
		In order to prove  the equality \ref{eq: toproove}  we need to show that 
		$$ \int_{\Omega_h} K_h\nabla \varphi_h \cdot (u_h \chi_h\nabla \eta_\delta) \to 0.$$
		We will use here that $F_{0,h}$ is the zero set of $u_h$.
		Note that $u_h\chi_h \in C^1_0(B(0,R_0))$ and vanishes on $F_{0,h}$, so $|u_h\chi_h| \leq A(u_h,\chi_h)\delta$ in the $\delta-$neighbourhood  of $F_{0,h}\cup C(0,R_0)$.  Thus $|u_h\chi_h\nabla \eta_{\delta}|$ is  bounded by $A(u_h,\chi_h)A_h$ in $B(0,R_0)$. Moreover $|(\textup{supp} \nabla \eta_{\delta})|=\mathcal{O}(\delta)$. Hence
		\begin{equation}
			\label{eq: divutil4}
			\int_{\Omega_h} K_h\nabla \varphi_h \cdot  (\chi_h u_h\nabla \eta_{\delta}) \leq A(u_h,\chi_h)A_h \sqrt{|(\textup{supp} \nabla \eta_{\delta})|} \sqrt{\int_{\Omega_h} |\nabla \varphi_h|^2} \to 0 \mbox{ as } \delta \to 0 .
		\end{equation}
		
		Thus with the equations \ref{eq: divutil2} ad \ref{eq: divutil3} and $\ref{eq: divutil4}$ we obtain
		$$\int_{\Omega_h} h^2K_h\nabla \varphi_h \cdot(\chi_h \nabla u_h  +u_h\nabla \chi_h) =  \int_{\Omega_h} K_hV_h\varphi_h u_h\chi_h=\int_{\Omega'_h} K_hV_h\varphi_h u_h\chi_h$$
		(the second equality is due to the fact that $u_h=0$ on $F_{0,h}$).
		
		With the extension of $\varphi_h$, we have $$\nabla \varphi_h=\nabla \tilde{K}_h^{-1/2}+ 1_{\Omega_h} \nabla \tilde{\varphi}_h \mbox{ on } B(0,R_0).$$
		So
		\begin{multline}
			\int_{\Omega_h} h^2K_h\nabla \varphi_h \cdot(\chi_h \nabla u_h  +u_h\nabla \chi_h)=\\
			\int_{\Omega_h} h^2K_h\chi_h\nabla \tilde{K}_h^{-1/2} \cdot \nabla u_h  +\int_{\Omega'_h}h^2K_h(\chi_h  1_{\Omega_h}  \nabla \tilde{\varphi}_h \cdot\nabla u_h  +u_h\nabla \varphi_h \cdot\nabla \chi_h),
		\end{multline}
		
		and, by integration by part, we can see that
		$$\int_{\Omega_h} h^2K_h\chi_h\nabla \tilde{K}_h^{-1/2} \cdot \nabla u_h=\int_{\Omega'_h} h^2K_h\chi_h\nabla \tilde{K}_h^{-1/2} \cdot \nabla u_h.$$

		Finally we get:
		$$\int_{\Omega'_h} h^2K_h\nabla \varphi_h \cdot(\chi_h \nabla u_h  +u_h\nabla \chi_h) =  \int_{\Omega'_h} V_h\varphi_h u_h\chi_h$$
		and, after subtracting the equation \ref{eq: divutil1}, we finish the proof of the inequality \ref{eq: toproove}.
		
	\end{proof}

	\subsection{Quasi-conformal mapping}
	\label{subsect: chg de var quasi-conforme}
	
	We dispose now of a function $f_h\in H_{loc}^{1}(\Omega'_h=B(0,R_0)\setminus F_{1,h})$ which is weak solution of the equation $-h^2\mathrm{div}(K_h\varphi_h^2\nabla f_h)=0$ in $\Omega'_h$.
	The purpose of this section is to show that up to a change of variable, $f_h$ is an harmonic function. For that we will use the theory of quasi-conformal mapping.
	
	First we use a stronger version of the Poincaré lemma:
	
	\begin{lem}
		\label{lem: Poincarelem}
		Let $F\in L_{loc}^p(B(0,1),\R^2)$ ($1\leq p<\infty$) a vector field such that $\mathrm{div}(F)=0$ in the sense of distribution. Then there is a function $u\in W_{loc}^{1,p}(B(0,1))$ such that $$F=(F_1,F_2):=(\partial_{x_2}u,-\partial_{x_1}u).$$
	\end{lem}
	The proof of this lemma can be find in Appendix.
	
	In our case, with this lemma \ref{lem: Poincarelem}, we dispose of a function $\tilde{f}_h\in H_{loc}^1$, define locally, such that
	$$\nabla\tilde{f}_h=\left( -K_h\varphi_h^2\,\partial_{x_2} f_h;  K_h\varphi_h^2\,\partial_{x_1} f_h \right).$$ Locally, $f_h$ and $\tilde{f}_h$ define a complex form $\omega_h:=f_h+i\tilde{f}_h$ which is solution of a Beltrami equation:
	
	\begin{equation}
		\label{eq: Beltrami}
		\frac{\partial \omega_h}{\partial \Bar{z}}=\mu_h \frac{\partial \omega_h}{\partial z}
	\end{equation}
	
	with coefficient:
	
	\begin{equation}
		\label{def: définition du cooef de beltrami}
		\mu_h := \frac{1-K_h\varphi_{h}^2}{1+K_h\varphi_{h}^2}.\frac{\partial_x f_{h} +i \partial_y f_{h}}{\partial_x f_{h} -i\partial_y f_{h}}.
	\end{equation}
	
	Unlike $\omega_h$, $\mu_h$ can be globally defined on $\Omega'_h$ and verifies:
	
	\begin{equation}
		|\mu_h|\leq \left|\frac{1-K_h\varphi_{h}^2}{1+K_h\varphi_{h}^2}\right|\lesssim \epsilon(h)\delta_{\|\nabla \phi\|_{\infty} \neq 0} +\mathcal{O}(\epsilon(h)^2).
	\end{equation}
	by construction of $\varphi_h$ (see proposition \ref{prop: constrsmallsolSchrod}).
	
	We extend $\mu_h$ by $0$ on the whole complex plan. Thus, $\mu_h$ is measurable function compactly supported with $|\mu_h|\lesssim  \epsilon(h)||\nabla \phi||_\infty +\mathcal{O}(\epsilon(h)^2(||\nabla \phi_h||_{\infty}^2 +||V_h-E_h||_{\infty}))<1$ for $\epsilon(h)$ small enough.  
	Consequently, we can apply an existence Theorem of quasi-conformal mapping to obtain solutions of Beltrami equations:
	
	\begin{prop}(\cite{ahlfors2006lectures}, chapter $5$, section $B$, Theorem $1$)
		\label{prop: def_de_psi_h}
		\hfill\\
		There is a quasi-conformal homeomorphism of the plan $\psi_h$, with dilatation coefficient  $T_h:=\frac{1+|\mu_h|}{1-|\mu_h|}$, such that:
		\begin{itemize}
			\item $\psi_h \in H_{loc}^{1}$,
			\vspace{0,2cm}
			\item $\frac{\partial \psi_h}{\partial \Bar{z}}=\mu_h \frac{\partial \psi_h}{\partial z}$,
			\vspace{0,2 cm}
			\item $T_h-1\lesssim \epsilon(h)\delta_{\|\nabla \phi\|_{\infty} \neq 0} +\mathcal{O}(\epsilon(h)^2)$,
			\vspace{0,2cm}
			\item $\psi_h(0)=0$ and $\psi_h(\infty)=\infty$.
			\vspace{0,2cm}
		\end{itemize}
	\end{prop}

	So we dispose of a quasi-conformal mapping $\psi_h$ of the plan onto itself which verify the same Beltrami equation as $\omega_h$.
	By a factorization Theorem of Stoilow (see Theorem $5.5.1$ of \cite{Astala2008-mi}), $\omega_h \circ \psi_h^{-1}$ is locally holomorphic and so harmonic.
	So $f_{h} \circ \psi_h^{-1}$ is in fact harmonic on $\psi_h(\Omega'_h)$ as the real part of $\omega_h \circ \psi_h^{-1}$. 
	
	By the Riemann uniformization Theorem, there exists a conformal biholomorphism mapping $\psi_h(B(0,R_0))$ on $B(0,R_0)$ and $\psi_h(0)$ on $0$.
	The composition of these two homeomorphisms is a $T_h-$quasi-conformal mapping $g_h$ satisfying the same Beltrami equation as $\omega_h$ and fixing $0$. Therefore $\tilde{u}_h:= f_h \circ g_h^{-1}$ is an harmonic function on $g_h(\Omega'_h)\subset B(0,R_0)$.

	\subsection{Space distortion}
	\label{subsect: distortion de lespace}
	In the previous section, we reduce the function $f_h$ to an harmonic one $\tilde{u}_h$ by sending the irregularities into the ambient metric. With the help of the quasi-conformal theory, we will see that the distortions are in fact well control.
	
	As a $T_h-$quasi-conformal mapping, the distances' distortions of $g_h$ inside $B(0,R_0)$ are Hölder:
	
	\begin{lem} (Mori Theorem, \citep{ahlfors2006lectures}, chapter $3$, section $C$)\\
		\label{lem: InegdeMori}
		for $z_1$,$z_2\in B(0,R_0)$.
		\begin{equation} 
			\frac{1}{16}\left| \frac{z_1-z_2}{R_0} \right|^{T_h} \leq \frac{|g_h(z_1)-g_h(z_2)|}{R_0}\leq16\left| \frac{z_1-z_2}{R_0} \right|^{1/T_h}
		\end{equation}	
	\end{lem}
	
	To be coherent with the different order of distances that we used previously (in particular in the section \ref{sect: decoupespace}), we would like that the images of distances proportional to $ h\epsilon$ to be also proportional $\simeq h\epsilon(h)$. For that, we know that $0\leq T_h-1\lesssim \epsilon(h)\delta_{\|\nabla \phi\|_{\infty} \neq 0} +\mathcal{O}(\epsilon(h)^2)$, thus we put:
	
	\begin{equation}
		\label{eq: epsilon(h)}
		\epsilon(h):=  \left\{
		\begin{array}{ll}
			\frac {c}{\log\left(R_0/h\right)} & \mbox{if } ||\nabla \phi||_\infty>0,\\
			\frac{c}{\log\left(R_0/h\right)^{1/2}} & \mbox{else.}
		\end{array}
		\right.
	\end{equation}
	
	Thus $(R_0)^{T_h} \sim R_0 \sim (R_0)^{1/T_h}$ and, if we take $c$ small enough, then if $h\epsilon(h)\leq |z_1-z_2|\leq 2R_0$ then
	
	\begin{equation}
		\label{eq: inegsurdistance}
		\frac{1}{32}\left| z_1-z_2 \right| \leq |g_h(z_1)-g_h(z_2)|\leq 32\left| z_1-z_2 \right|.
	\end{equation}
	It should be noted that with the choice of $\epsilon(h)$, $c$ can be chosen independently of $h$.
	
	Hence, the function $\tilde{u}_h:=f_{h} \circ g_h^{-1}$ is harmonic on $B(0,R_0)\setminus  g_h(F_{1,h})$ where $g_h(F_{1,h})$ is disjoint union of sets with diameter proportional to $h\epsilon(h)$. Each set of $g_h(F_{1,h})$ is contained in a ball with radius $32h\epsilon(h)$. Let's note these balls $\tilde{D}_{j,h}$. Each set of $g_h(F_{1,h})$ is separated by at least a distance of $\frac{C}{32}h\epsilon(h)$ from  the others sets. So $\tilde{D}_{j,h}$ are separated between themselves, from $C(0,R_0)$ and from $\{\tilde{f}_h=0\}$  by a distance bigger than $\frac{C}{32}h\epsilon(h)-128h\epsilon(h)=:32\tilde{C}h\epsilon(h)$, with:
	
	\begin{equation}
		\tilde{C}:=\frac{C}{32^2}-4.
	\end{equation}
	For the following, we will suppose that $\tilde{C}>100$ (which is true for $C$ large enough).

	\section{End of Theorem \ref{thm: mainthm}'s proof}
	\label{sect: end_proof_main}
	
	In the last section, we have substitute, as done in \cite{logunov2020landisconjectureexponentialdecay}, the eigenfunction $u_h$ by an harmonic function $\tilde{u}_h$. So we've reached the same point as they have. This means we'll be able to rely on their "local main Theorem": 
	
	\begin{thm} (\cite{logunov2020landisconjectureexponentialdecay}, Theorem $5.3$.)\\
		\label{thm: Toy}
		Let $D_j$ be a collection of $100-$separated unit disks on $\mathbb{R}^2=\mathbb{C}$ such that 
		$0 \notin \cup 3D_j$. Let $R> 10^4$, $0<r\leq R/4$. 
		Consider any harmonic function $u$ in $B(0,R) \setminus \cup D_j$ such that $u$ does not change sign in
		$(5D_j\setminus D_j) \cap B(0,R)$ for every $j$. Assume  that 
		$$\sup \limits _{B(0,R-R/64) \setminus \cup 3D_j}|u| \geq e^{- N} \sup \limits_{B(0,R) \setminus \cup 3D_j}|u|.$$
		Then 
		\begin{equation}
			\sup \limits _{B(0,r) \setminus \cup 3D_j}|u| \geq  \left(\frac r R \right)^{ L (R + N) } \sup \limits_{B(0,R) \setminus \cup 3D_j}|u|,
		\end{equation}
		with some absolute constant $L>0$.
	\end{thm}
	
	For completeness, this Theorem \ref{thm: Toy} is proved in Section \ref{sect: thm_Toy_proof}.
	
	In order to apply this Theorem, we will rescale the distances by the factor $(32h\epsilon(h))^{-1}$. Thus $R_h=R_0(32h\epsilon(h))^{-1}$ and $D_{j}:= (32h\epsilon(h))^{-1}\tilde{D}_{j,h}$ are a collection of unit disk $\tilde{C}-$separated between themselves and from the zero-set of $\tilde{u}_h(\cdot (32h\epsilon(h))^{-1})$. Therefore, we deduce from the Theorem \ref{thm: Toy}:
	
	If $$ \sup \limits _{B(0,R_0-R_0/64) \setminus \cup 3\tilde{D}_{j,h}}|\tilde{u}_h| \geq e^{- \tilde{N}_h} \sup \limits_{B(0,R_0) \setminus \cup 3\tilde{D}_{j,h}}|\tilde{u}_h|,$$ then
	\begin{equation}
		\label{eq: appthLog1}
		\sup \limits _{B(0,\tilde{r})\setminus \cup 3\tilde{D}_{j,h}}|\tilde{u}_h| \geq  \left(\frac{\tilde{r}}{R_0} \right)^{ L (R_h + \tilde{N}_h) } \sup \limits_{B(0,R_0) \setminus \cup 3\tilde{D}_{j,h}}|\tilde{u}_h|,
	\end{equation}
	
	where $\tilde{r}<R_0/4$.
	
	In the section \ref{sect: decoupespace}, we specifically asked that points $0$ and $x_{\max,h}$ (the point where $\sup_{B(0,R_0/2)}|u_h|$ is attained) to be $Ch\epsilon(h)-$separated from the disks of $F_{1,h}$. By construction of the quasi-conformal homeomorphism, $g_h(0)=0$, so the disks $\tilde{C}\tilde{D}_{j,h}$ do not contain $0$ and $g_h(x_{\max,h})$. The distortion estimate implies that $g_h(x_{\max,h}) \in \overline{B(0,R_0-R_0/64)}$. Moreover, by the construction of $\varphi_h$ (proposition \ref{prop: constrsmallsolSchrod}) and by the equations \ref{eq: convsmallsol1} and \ref{eq: convsmallsol2}, 
	\begin{equation}
		\label{eq: contrlfinalsurlafunharm}
		\frac{(\sqrt{K_h}\,|u_h|) }{2}\circ g_h^{-1}\leq |\tilde{u_h}|\leq 2(\sqrt{K_h}\,|u_h|) \circ g_h^{-1}
	\end{equation}
	 for $h$ small enough. Therefore, if 
	$$\frac{\sup_{B(0,R_0)}|u_h|}{\sup_{B(0,R_0/2)}|u_h|} \leq e^{N_h},$$
	then
	$$ \frac{\sup_{B(0,R_0) \setminus \cup 3\tilde{D}_{j,h}}|\tilde{u}_h|}{\sup_{B(0,R_0-R_0/64)\setminus \cup 3\tilde{D}_{j,h}}|\tilde{u}_h|} \leq 4e^{\mathrm{osc}(\phi_h)/2h}e^{N_h}=: e^{\tilde{N}_h},$$
	where $\mathrm{osc}(\phi_h):=\max \phi_h-\min \phi_h$.
	
	For a ball $B(0,r)$ in $B(0,R_0)$, we know by the Mori lemma that $g_h(B(0,r))$ contains  a disk with centre at $0$ and radius $$\frac{R_0}{16}\left(\frac{r}{R_0}\right)^{T_h} \geq \frac{R_0}{16} \left(\frac{r}{R_0}\right)^2=:\tilde{r}.$$
	This imply that for $r\leq R_0/4$, we have $\tilde{r}\leq R_0/4$. Therefore, from the equation \ref{eq: appthLog1}, we get
	
	\begin{equation}
		\label{eq: respartielharmcontrlgrandparpetit}
		\sup \limits _{B(0,\tilde{r})} |\tilde{u}_h| \geq  \left(\frac{r}{4R_0} \right)^{ 2L (R_h \tilde{N}_h) } \sup \limits_{B(0,R_0-R_0/64) \setminus \cup 3D_{j,h}}|\tilde{u}_h|.
	\end{equation}
	
	Now, by the fact $g_{h}(x_{max})\in B(0,R_0-R_0/64) \setminus \cup 3\tilde{D}_{j,h}$ and by the inequalities \ref{eq: contrlfinalsurlafunharm} : 
	\begin{equation}
		\label{eq: boudbybelowharmbymaxofschrod}
			\sup \limits_{B(0,R_0-R_0/64) \setminus \cup 3\tilde{D}_{j,h}}|\tilde{u}_h|\geq  \frac{e^{\min(\phi_h)/2h}}{2}\,\sup \limits_{B(0,R_0/2)}|u_h|\geq \frac{e^{\min(\phi_h)/2h}}{2}e^{-N_h}\,\sup \limits_{B(0,R_0)}|u_h|.
	\end{equation}
 Thus by the equations \ref{eq: respartielharmcontrlgrandparpetit} and \ref{eq: boudbybelowharmbymaxofschrod}  we get:
	
	\begin{equation}
		\sup \limits _{B(0,r)} |u_h| \geq \frac{1}{4}e^{-\mathrm{osc}(\phi_h)/2h}\left(\frac{r}{4R_0} \right)^{ 2L (R_h + \tilde{N}_h) } e^{-N_h}\,\sup \limits_{B(0,R_0)}|u_h|.
	\end{equation}
	
	Thus for $h$ small enough and because $r/R_0\leq 4^{-1}\leq e^{-1}$, we have
	\begin{equation}
		\sup \limits _{B(0,r)} |u_h| \geq \left(\frac{r}{R_0} \right)^{ L'(R_0 (h\epsilon(h))^{-1} + N_h+\mathrm{osc}(\phi_h)/2h) } \,\sup \limits_{B(0,R_0)}|u_h|
	\end{equation}
	with $L'>0$ absolute constant. Eventually, we obtain the main Theorem \ref{thm: mainthm}.
	
	\section{Localization on compact setting.}
	
	Now that we have proved the Theorem  \ref{thm: mainthm}, we can study the eigenfunctions of weighted Schrödinger operator $P_{h}$ on compact Riemann surface. Let be $(M,g)$ a $W^{2,\infty}$ compact connected Riemann surface equipped with a Lipschitz (definite positive) metric. Let consider the weighted Schrödinger eigenvalue equation:
	
	\begin{equation}
		\label{eq: Schrodinger equation_ClManifold}
		P_hu_h:= \left(-\frac{h^2}{K_h}\mathrm{div}_{g}(K_h\nabla_g) +V_h\right)u_h=E_hu_h \quad \mbox{ on } M,
	\end{equation}
	
	where 
	\begin{itemize}
		\item $K_h:=\exp(\phi_h/h) \mbox{ with } \phi_h=\phi+\mathcal{O}_{W^{1,\infty}}(h) \in W^{1,\infty}(M,\R),$
		\item $V_h:=V+\mathcal{O}_{L^\infty}(h)\in L^\infty(M,\R)$,
		\item and $E_h:=E+\mathcal{O}_{\R}(h)\in \R$,
	\end{itemize}
	with $h$ the semiclassical parameter.
	
	\subsection{Closed Surface}
	\label{sect: closed_surface}
	
	Assume in this section that the manifold $M$ is closed.
	 
	\begin{thm}
		\label{thm: closed_Riem_surf}
		For $h$ small enough, let be  $u_h\in H^1(M,\R)$ weak solution of equation \ref{eq: Schrodinger equation_ClManifold}.
		Then for any  ball $B_r$ of radius $r<1/4$:
		
		\begin{equation}
			\sup\limits_{B_r}|u_h| \geq  r^{ L\left(\frac{\epsilon(h)+\mathrm{osc}(\phi_h)}{h}\right)} \sup_{M}|u_h|,
		\end{equation}
		where $L$ is an absolute constant, $\mathrm{osc}(\phi_h)= \max_M \phi_h -\min_M \phi_h$  and 
		$$
		\epsilon(h):=  \left\{
		\begin{array}{ll}
			\Theta\log\left(\frac{\Theta}{h}\right) & \mbox{if } 	||\nabla \phi||_{\infty}>0,\\
			\Theta\log\left(\frac{\Theta}{h}\right)^{1/2} & \mbox{else.}
		\end{array}
		\right.
		$$ 
		with $\Theta= \max \{\|\nabla \phi\|_{\infty}, \sqrt{\|V-E\|_{\infty}},1\}$. 
	\end{thm}
	
	\begin{proof}
		In fact, specifically in dimension $2$, we can use isothermal coordinates to rewrite the equation $$-\frac{h^2}{K_h}\mathrm{div}_g( K_h\nabla_g )u_h + V_h u_h =E_hu_h$$ as $$-\frac{h^2}{\tilde{K}_h}\mathrm{div}( \tilde{K}_h\nabla )\tilde{u}_h + \tilde{V}_h \tilde{u}_h =0$$ with the euclidean metric.
		 Indeed, let be $(\Psi_j:U_j \to \R^2)_j$ an atlas of $W^{2,\infty}-$diffeomorphisms (in the sense that the change of charts are in $W^{2,\infty}\cap \mathcal{C}^1$)  and  let be $I_g $ the first fundamental form associated to the metric $g$. In local coordinates $(x,y)\in \R^2$, $I_g$ take the form $$I_g=g_{1,1}dx^2+2g_{1,2}dxdy+g_{2,2}dy^2$$ with $g_{i,j}\in W^{1,\infty}(\R^2)$. In isothermal coordinate $(x',y')\in \R^2$, $I_g$ should look like $$I_g=e^\psi (dx'^2+dy'^2)$$ with $\psi$ continuous. 
		Let switch to complex coordinates $z=x+iy$ and $z'=x'+iy'$. we note that
		$$ I_g=\lambda |dz+\mu d\overline{z}|^2 $$
		with $$\lambda=\frac{1}{4}(g_{1,1}+g_{2,2}+2\sqrt{g_{1,1}g_{2,2}-g_{1,2}^2}); \quad \mu= \frac{g_{1,1}-g_{2,2}+2ig_{1,2}}{4\lambda}$$
		and $$e^{\psi}|dz'|^2=e^{\psi}|\partial_z z'|^2\left|dz+ \frac{\partial_{\overline{z}}z'}{\partial_z z'}d\overline{z}\right|.$$
		
		 Thus our isothermal changes of variables are obtain locally as solution of the Beltrami equation
		 \begin{equation}
		 	\label{eq: beltamiiso}
		 	\partial_{\overline{z}} z'=\mu\, \partial_z z'
		 \end{equation}
		 and $e^{\psi} |\partial_z z'|^2=\lambda$.
		 
		 The Beltrami coefficient $\mu$ verify $\mu \in W^{1,\infty}(\C)$ because $g_{i,j}\in W^{1,\infty}(\R^2)$ and $\sqrt{g_{1,1}g_{2,2}-g_{1,2}^2}$ is positively bounded by below, and $|\mu|\leq 1$ on $\C$. By the Theorem $1$ of \cite{ahlfors2006lectures} (chapter $5$, section $B$), we dispose of an homeomorphism $\Theta$ of $\C$ into itself, solution of the Beltrami equation \ref{eq: beltamiiso}. Moreover by the lemma $3$ of \cite{ahlfors2006lectures} (chapter $5$, section $B$), the mapping $\Theta: z\in \C\to z' \in\C$ is $\mathcal{C}^1$ and so $\psi$ is continuous. Furthermore, by the factorization Theorem of Stoilow (see Theorem $5.5.1$ of \cite{Astala2008-mi}), two solutions of the same Beltrami equation are the same up to a composition by an holomorphic function.  Thus the charts $\Theta_j \circ\Psi_j$ verify that their change of charts  are still in $W^{2,\infty}\cap \mathcal{C}^1$. Thus $M$ conserve its $W^{2,\infty}$ structure.
		 
		  So after isothermal change of variable, we reduce locally to $$-\frac{h^2}{\tilde{K}_h}\mathrm{div}( \tilde{K}_h\nabla )u_h + e^{\psi}(\tilde{V}_h-E_h) u_h =0$$ with $\tilde{K}_h\in W^{1,\infty}$, $\tilde{V}_h\in L^\infty$ and $e^\psi\in \mathcal{C}^0$.
		Now, let $x_{\max}$ be the point in $M$ where $|u_h|(x_{\max})=\sup_M |u_h|$ and take $x\in M$. The idea of the proof is to link $x_{\max}$ to $x$ by a succession of disks with radius $1$. The number of disks needed is controlled by the compacity of $M$. We start with a ball centre on $x_1\neq x_0=x_{\max}$ such that $|x_1-x_0|\leq 1/2$ so $$\sup_{B(x_1,1)} |u_h|=   \sup_{B(x_1,1/2)} |u_h|.$$ Thus by the Theorem \ref{thm: mainthm}
		
		\begin{equation}
			\sup \limits _{B(x_1,1/4)} |u_h| \geq \left(\frac{1}{4} \right)^{L_1\left(\frac{\epsilon(h)+\mathrm{osc}(\phi_h)}{h}\right) }\,\sup \limits_{M}|u_h|.
		\end{equation}
		
		Now if you take $x_2$ such that $B(x_1,1/4) \subset B(x_2,1/2)$, you will have$$\frac{\sup_{B(x_2,1)} |u_h|}{\sup_{B(x_2,1/2)} |u_h|}\leq \frac{\sup_{M} |u_h|}{\sup_{B(x_1,1/4)} |u_h|}\leq e^{ \log(4)L_1\left(\frac{\epsilon(h)+\mathrm{osc}(\phi_h)}{h}\right)}.$$ Hence
		
		\begin{align}
			&\sup \limits _{B(x_2,1/4)} |u_h| \geq \left(\frac{1}{4} \right)^{L_1\left(\left(\frac{\epsilon(h)+\mathrm{osc}(\phi_h)}{h}\right)+\log(4)L_1\left(\frac{\epsilon(h)+\mathrm{osc}(\phi_h)}{h}\right)\right)}\,\sup \limits_{B(x_2,1)}|u_h|\\
			&\geq \left(\frac{1}{4} \right)^{L_2\left(\frac{\epsilon(h)+\mathrm{osc}(\phi_h)}{h}\right)}\,\sup \limits_{B(x_1,1/4)}|u_h|\\
			&\geq  \left(\frac{1}{4} \right)^{(L_2+L_1)\left(\frac{\epsilon(h)+\mathrm{osc}(\phi_h)}{h}\right)}\,\sup \limits_{M}|u_h|.
		\end{align}
		
		And so we can repeat the process until $x_n\in B(x,1/2)$, to get:
		\begin{equation}
			\sup \limits _{B(x,r)} |u_h| \geq r^{\left(\sum_{k=1}^{n} L_k\right)\left(\frac{\epsilon(h)+\mathrm{osc}(\phi_h)}{h}\right)}\,\sup \limits_{M}|u_h|.
		\end{equation}
	\end{proof}
	
	\subsection{Compact surface with boundary conditions}
	\label{sect: compact_surface_with_boundary}
	
	In this section, we show that the precedent result can be extended to compact surface with boundary on which we put Robin conditions:
	\begin{thm}
		For $h$ small enough, let be  $u_h\in H^1(M,\R)$ weak solution of the equation \ref{eq: Schrodinger equation_ClManifold}. Moreover, assume $u_h$ verify the Robin boundary condition: $$h\,\alpha(x)\,\partial_n {u_h}_{|\partial M}(x)+\beta(x)\,{u_h}_{|\partial_M}(x)=0$$ with $n$ the normal exterior to the boundary and $\alpha,\,\beta\in W^{2,\infty}(\partial M, \R)$ such that 
		\begin{itemize}
			\item if $\alpha(x)=0$ for $x\in \partial M$, then $\alpha=0$ on the connected component of $\partial M$ that contained $x$,
			\vspace{1mm}
			\item $\frac{\beta}{\alpha}\in W^{1,\infty}(\mathrm{Supp}(\alpha))$,
			\vspace{1mm}
			\item $\{\alpha=0\}\cap\{\beta=0\}=\emptyset$.
		\end{itemize}
		
		 Then, for any  ball $B_r\subset \overline{M}$ of radius $r<1/4$:
		
			\begin{equation}
			\sup\limits_{B_r}|u_h| \geq  r^{ L\left(\frac{\epsilon(h)+\mathrm{osc}(\phi_h)+\sup\left|\frac{\beta}{\alpha}\right|1_{\{\alpha \neq 0\}}}{h}\right)} \sup_{M}|u_h|
		\end{equation}
			where $L$ is an absolute constant, $\mathrm{osc}(\phi_h)= \max_M \phi_h -\min_M \phi_h$  and 
		$$
		\epsilon(h):=  \left\{
		\begin{array}{ll}
			\Theta\log\left(\frac{\Theta}{h}\right) & \mbox{if } 	||\nabla \phi||_{\infty}+\| \frac{\beta}{\alpha} 1_{\{\alpha \neq 0\}}\|_{W^{1,\infty}}>0,\\
			\Theta\log\left(\frac{\Theta}{h}\right)^{1/2} & \mbox{else.}
		\end{array}
		\right.
		$$ 
		with $\Theta= \max \{\|\nabla \phi\|_{\infty}+\| \frac{\beta}{\alpha} 1_{\{\alpha \neq 0\}}\|_{W^{1,\infty}}, \sqrt{\|V-E\|_{\infty}},1\}$. 
	\end{thm}
	
	We want to emphasize that Robin boundary conditions play an exponential role in our result. However, Dirichlet and Neumann boundary conditions, which are types of Robin boundary conditions, allow us to avoid this form of exponential weight. In particular, in the context of the Schrödinger operator with such boundary conditions, we obtain the following subresult:
	
	\begin{cor}
		Let $(M,g)$ be a compact Riemannian surface with a $W^{2,\infty}$ structure with $g\in W^{1,\infty}(M)$. For $h$ small enough, let be  $u_h\in H^1(M,\R)$ weak solution of $$-h^2\Delta_g u_h + V_h u_h =E_hu_h \mbox{ on } M$$ where  $V_h=V+\mathcal{O}_{L^\infty}(h)\in L^\infty(M,\R)$ and $E_h= E +\mathcal{O}_{\R}(h) \in\R$. Moreover, assume $u_h$ verify the boundary condition: $$ {u_{h}}_{|\partial M}=0 \mbox{ (Dirichlet)}\quad \mbox{ or }\quad {\partial_nu_{h}}_{|\partial M}=0 \mbox{ (Neumann)}$$ with $n$ the normal exterior to the boundary. 
		Then for any  ball $B_r\in \overline{M}$ of radius $r<1/4$:
		
		\begin{equation}
			\sup\limits_{B_r}|u_h| \geq  r^{ \frac{L\max \{\sqrt{\|V-E\|_{\infty}},1\}}{h}\log\left(\frac{\max \{\sqrt{\|V-E\|_{\infty}},1\}}{h}\right)^{1/2}} \sup_{M}|u_h|
		\end{equation}
		where $L$ is an absolute constant.
	\end{cor}
	
	\begin{proof}
		 To prove this last Theorem, we will try to reduce to closed case by gluing the manifold on copies of itself along their boundaries in a way to conserve the Lipschitz regularity. However,  the usual method for this procedure relies on the use of a reflexion principle in geodesic coordinate systems. At our level of regularity, the geodesic coordinate system doesn't exist. To overcome this difficulty, we will use the pseudo-geodesic system define by \cite{Burq2022}. This system relies on a regularisation of the normal direction to the boundary (see the proof of the following lemma in Appendix).
		 
		\begin{lem}
			\label{lem: local_pseudogeod_coord}
			Assume that $M$ is a $W^{2,\infty}$ Riemannian surface with boundary (the regularity here means that changes of charts are in $W^{2,\infty}\cap \mathcal{C}^1$), endowed with a Lipschitz (definite positive) metric $g$ and a Lipschitz (positive) density $K$. Let
			$$ \Delta = \frac{1}{K}\mathrm{div}(g^{-1}K\nabla).$$
			Then near any point of the boundary $x_0\in \partial M$ there exists a $W^{2,\infty}$ coordinate system $\Psi_0:U_0\subset \overline{M} \to B(0,\delta)\cap (\R_x \times[0, +\infty))$ with $\delta>0$ such that
			in this coordinate system $x_0 = (0, 0) $ and $$\Delta=\frac{1}{K}\nabla^{T}K\,b(x,y)\nabla,\quad b_{|\partial \Omega}=\begin{pmatrix}
				b_0(x) & 0\\
				0 & 1
			\end{pmatrix}$$
			with $b_0\in W^{1,\infty}(\partial\Omega)$ positive. 
		\end{lem}

		Let $M'=\overline{M}\times \{\pm 1\}/ \{\alpha=0\}$ (or $M$ if $\{\alpha=0\}$ is empty) and $\tilde{M}=\overline{M'}\times \{\pm i\}/\partial M'$ (or $M'$ if $\partial M'$ is empty). With this lemma \ref{lem: local_pseudogeod_coord}, we can define a covering of $\partial M \subset \bigcup_{j=1}^N U_j$ (seen as a subset of $\overline{M}$) and $W^{2,\infty}$ diffeomorphisms $\Psi_j: U_j \to \R^2=\R_x\times \R_y$ 
		such that $\Psi_j(U_j)\subset B(0,1)\cap \{y\geq 0\}$ and the metric on $\{\alpha=0\}$  takes the form 
		\begin{equation}
			\label{eq: form_metric_boundary}
			b'_{|\partial \Omega}(x,y)=\begin{pmatrix}
				b'(x) & 0\\
				0 & 1
			\end{pmatrix}
		\end{equation}
		  with $\tilde{b}\in W^{1,\infty}$.
		
		We can now perform the first gluying by defining a covering of $\{\alpha=0\}$ (now seen as a subset of $M'$ ):
		$$\{\alpha=0\}\subset \bigcup_{j=1}^{N}U_j\times\{\pm 1\}=\bigcup_{j=1}^{N}U'_j$$
		where we identify the points in $\left(U_j\cap \{\alpha=0\}\right)\times \{\pm 1\}$, and define the map
		$$\Psi'_j:z=(x,\epsilon)\in U_j\times \{\pm 1\}\to \left\{
			\begin{array}{ll}
					\Psi_j(x) & \mbox{ if } \epsilon=1, \\
				S\circ \Psi_j(x) & \mbox{ if } \epsilon=-1,
			\end{array}
			\right.$$
			where $S(x,y)=S(x,-y)$.
			We conclude that the image of the metric induced on $M'$ by the metrics on the two copies of $M$ is well defined and Lipschitz. Indeed, on $\Psi'_j(U_j\times\{1\})$ the metric is given by $$b_{+1}(x,y)1_{\{y\geq0\}},$$
			while on $\Psi'_j(U_j\times\{-1\})$ the metric is given by
			$$
			b_{-1}(x,y)= \begin{pmatrix}
				1 &0\\
				0 &-1
			\end{pmatrix}b \circ S(x, y)\begin{pmatrix}
			1 &0\\
			0 &-1
			\end{pmatrix}$$
			As a consequence, it follows from \ref{eq: form_metric_boundary} that the two metrics coincide on $\{y=0\}$  and they define a Lipschitz metric $b'$ on a small neighbourhood of $(0,0)\in (-\delta,\delta)_x\times (-\epsilon,\epsilon)_y$ with $\epsilon,\, \delta>0$. For the $W^{2,\infty}$ regularity, let derive the changing of chart: for $j\neq k$, 
			$$\mathrm{d}_{x,y} \Psi'_{j} \circ (\Psi'_{k})^{-1}=\left\{\begin{array}{ll}
				\mathrm{d}_{x,y}\Psi_j \circ \Psi_k^{-1} & \mbox{ on } y\geq 0,\\
				\begin{pmatrix}
					1 &0\\
					0 &-1
				\end{pmatrix} \mathrm{d}_{x,y}(\Psi_j \circ \Psi_k^{-1}) \begin{pmatrix}
				1 &0\\
				0 &-1
				\end{pmatrix} & \mbox{ on } y\leq 0,
			\end{array}\right. $$
		where $\mathrm{d}_{x,y}$ is the differential. By definition of $\Psi_j$ and construction of the metric $b'$ on the boundary, $\mathrm{d}_{x,y} \Psi'_{j} \circ (\Psi'_{k})^{-1}$ send the normal vector 
		$\left(\begin{array}{ll}
			0\\
			1
		\end{array}\right)$
		to the normal vector 
		$\left(\begin{array}{ll}
			0\\
			1
		\end{array}\right)$
		and tangents vectors to tangents vectors,  so on $\{y=0\}$: 
		\begin{equation}
			\label{eq: form_changechart_boundary}
			\mathrm{d}_{x,y}{\Psi_j \circ \Psi_k^{-1}}_{|\{y=0\}}(x,0)=\begin{pmatrix}
				q(x) & 0\\
				0 & 1
			\end{pmatrix}
		\end{equation}
		with $q\in W^{1,\infty}\cap \mathcal{C}^0$.
		We deduce that the two limits of the differentials $y\to 0^{\pm}$ coincide and consequently, the differential is $\mathcal{C}^0$. Let us now study the $L^\infty$ boundedness of derivatives of second order. The derivatives $\mathrm{d}_{x,x}^2$ and $\mathrm{d}_{x,y}^2$ follow from  \ref{eq: form_changechart_boundary} by taking a $\mathrm{d}_x$ derivative. Finally, the case of $\mathrm{d}_{y,y}^2$ derivative follows from the jump formula and the use of \ref{eq: form_changechart_boundary} which has shown
		that the first order derivatives have no jump, as
		$$\mathrm{d}_{y,y}^2\Psi'_{j} \circ (\Psi'_{k})^{-1}= \frac{\partial^2}{\partial_x^2}\left(\Psi_{j} \circ (\Psi_{k})^{-1}\,1_{\{y>0\}}+S\circ \Psi_{j} \circ (\Psi_{k})^{-1}\circ S\,1_{y<0}\right)$$
		$$=\frac{\partial^2}{\partial_x^2}\left(\Psi_{j} \circ (\Psi_{k})^{-1}\right)1_{\{y>0\}}+\frac{\partial^2}{\partial_x^2}\left(S\circ \Psi_{j} \circ (\Psi_{k})^{-1}\circ S\right)1_{y<0}$$
		$$+2\left(\frac{\partial}{\partial_x}\left(\Psi_{j} \circ (\Psi_{k})^{-1}\right)-\frac{\partial}{\partial_x}\left(S\circ \Psi_{j} \circ (\Psi_{k})^{-1}\circ S\right)\right)_{|y=0}\otimes\delta_{y=0}$$
		$$+\left(\Psi_{j} \circ (\Psi_{k})^{-1}-S\circ \Psi_{j} \circ (\Psi_{k})^{-1}\circ S\right)_{|y=0}\otimes\delta'_{y=0}$$
		$$=\frac{\partial^2}{\partial_x^2}\left(\Psi_{j} \circ (\Psi_{k})^{-1}\right)1_{\{y>0\}}+\frac{\partial^2}{\partial_x^2}\left(S\circ \Psi_{j} \circ (\Psi_{k})^{-1}\circ S\right)1_{y<0}.$$
		Hence the $W^{2,\infty}$ regularity.
		
		We dispose now of a $W^{2,\infty}$ structure on $M'$ with a Lipschitz metric $g'$ such that it coincide with $g$ away from $\{\alpha=0\}$. We extend the function $u_h$ on $M'$ by 
		\begin{equation}
			u'_h(z=(x,\pm1))=\pm u_h(x).
		\end{equation}
		This definition makes sense because $u_h(\{\alpha=0\})=0$. we also extend $V_h$ by $V'_h(x,\pm1)=V_h(x)$ and $\phi_h$ by $\phi'_h(x,\pm1)=\phi_h(x)$ which is in $W^{1,\infty}$ by the jump formula. Now we check that $u'_h$ verify
		$$P'_h u'_h:=\left(-\frac{h^2}{K'_h}\mathrm{div}_{g'}(K'_h\nabla_{g'})+V'_h\right)u'_h= E_h u'_h \mbox{ inside } M'.$$ Away from   $\{\alpha=0\}$ this is clear while near a point $x\in \{\alpha=0\}$ we can work in a coordinate chart $(\Psi'_j, U'_j)$. In this coordinate chart, the function
		$u'_h$ is defined by
		$$ u'_h=\left\{
			\begin{array}{ll}
			u_h(x,y) &\mbox{ if } y>0,\\
			-u_h(x,-y) &\mbox{ if } y<0.
		\end{array}
		\right.$$
		In $\{\pm y>0\}$, $u'_h$ satisfies the equation, and near $\{\alpha=0\}$ in our coordinate systems, we have
		$$u'_h(x,y)=u_h(x,y)\,1_{y>0}\,-u_h(x,-y)\, 1_{y<0},$$ and we have
		$$\partial_x u'_h(x,y)=\partial_x (u_h(x,y))\,1_{y>0}\,-\partial_x(u_h(x,-y))\, 1_{y<0};$$
		$$\partial_y u'_h(x,y)=\partial_y (u_h(x,y))\,1_{y>0}\,-\partial_y(u_h(x,-y))\, 1_{y<0}\,+(u_h(x,0^{+})+u_h(x,0^{-}))\otimes \delta_{y=0}$$
		$$=\partial_y (u_h(x,y))\,1_{y>0}\,-\partial_y(u_h(x,-y))\, 1_{y<0},$$
		because $u_h(x,0^{\pm})=0$ for all $x$. Consequently, the extension define a function in $H^1(M',\R)$. Then with the form of $g'$ on $\{\alpha=0\}$ and by the jump formula:
		\begin{align}
			h^2{K'_h}^{-1}\mathrm{div}_{g'}(K'_h\nabla_{g'})u'_h(x,y)&= (V'_h(x,y)-E_h)u'_h(x,y)\nn\\
			&+\left(\partial_y u_h(x,0^{+})-\partial_y u_h(x,0^{-})+0\times (\partial_x u_h(x,0^{+})+\partial_x u_h(x,0^{-})) \right)\otimes \delta_{y=0}\nn\\
			&=(V'_h(x,y)-E_h)u'_h(x,y).
		\end{align}
		
		So $u'_h\in H^1(M',\R)$ is solution of $P'_h u'_h= E_h u'_h \mbox{ inside } M'$. Moreover on $\partial M'$, $u'_h$ is solution of $$h\partial_n {u'_h}_{|\partial M'}+\gamma \, {u'_h}_{|\partial M'}=0$$ with
		\begin{equation}
		\gamma(z=(x,\pm 1))=\frac{\beta}{\alpha}(x)\in W^{2,\infty}(\partial M',\R).
		\end{equation}
		 The function $\gamma$ is well defined on $\partial M'$ by the hypothesis on $\alpha$ and $\beta$. 
		
		Let $\mathcal{U}$ be a tubular neighbourhood of $\partial M'$ in $M'$. In normal coordinate $(x,r)\in \mathcal{U}$ we define 
		\begin{equation}
			\mu(x,r):=\chi(r)\,re^{-\langle r\, \mathcal{D}_x\rangle+1}\gamma(x)
		\end{equation}
		where $\chi$ is a smooth nonnegative function on $U$ equal to $1$ on a $\epsilon-$neighbourhood of $\partial M'$ and is zero outside a $2\epsilon-$neighbourhood of $\partial M'$. We extend $\mu$ to $M'$ by zero.
		
		\begin{lem}
			The function $\mu$  is in $W^{2,\infty}(M',\R)$ and $\partial_{n}\mu_{|\partial M'}=\gamma$.
		\end{lem}
		\begin{proof}
			With the help of lemma \ref{lem: regularization_protocol}, we see that $\mu\in W^{2,\infty}(M',\R)$. Indeed,
			$$\partial_x \mu(x,r)=\chi(r)\,re^{-\langle r\, \mathcal{D}_x\rangle+1}(\partial_x\gamma)(x)\in L^\infty$$
			$$\partial_r \mu (x,t)=(\partial_r\chi)(r)\,re^{-\langle r\, \mathcal{D}_x\rangle+1}\gamma(x)+\chi(r)\,e^{-\langle r\, \mathcal{D}_x\rangle+1}\gamma(x)-\chi(r)\,r^2\mathcal{D}_x^2\langle r\, \mathcal{D}_x\rangle^{-1} e^{-\langle r\, \mathcal{D}_x\rangle+1}\gamma(x)$$
			$$=(r(\partial_r \chi)(r) +\chi(r)(1-\langle r\, \mathcal{D}_x\rangle+\langle r\, \mathcal{D}_x\rangle^{-1}))e^{-\langle r\, \mathcal{D}_x\rangle+1}\gamma(x)\in L^\infty.$$
			Moreover, we see that
			$$\partial_r \mu_{|\partial_M'}=\gamma.$$
			Then  for second derivatives:
			$$\partial_x^2 \mu(x,r)=\chi(r)\,(r\partial_{x})e^{-\langle r\, \mathcal{D}_x\rangle+1}(\partial_{x}\gamma)(x)\in L^\infty,$$
			$$\partial_x\partial_{r}\mu(x,r)=(r(\partial_r \chi)(r) +\chi(r)(1-\langle r\, \mathcal{D}_x\rangle+\langle r\, \mathcal{D}_x\rangle^{-1}))e^{-\langle r\, \mathcal{D}_x\rangle+1}\partial_{x}\gamma(x)\in L^\infty,$$
			$$\partial_r^2 \mu(x,r)=r\mathcal{D}_x \chi(r)(1-2\langle r\, \mathcal{D}_x\rangle^{-1}-\langle r\, \mathcal{D}_x\rangle^{-2}-\langle r\, \mathcal{D}_x\rangle^{-3})e^{-\langle r\, \mathcal{D}_x\rangle+1}(\partial_{x}\gamma)(x)$$
			$$+(\partial_r^2\chi)(r)\,re^{-\langle r\, \mathcal{D}_x\rangle+1}\gamma(x)+2(\partial_r\chi)(r)(1-\langle r\, \mathcal{D}_x\rangle+\langle r\, \mathcal{D}_x\rangle^{-1})e^{-\langle r\, \mathcal{D}_x\rangle+1}\gamma(x)\in L^\infty.$$
			Thus $\mu \in W^{2,\infty}$.
		\end{proof}
		
		Note that with the proof and the lemma \ref{lem: regularization_protocol},  we have that $\mathrm{osc}(\mu)\lesssim 2\| \gamma \|_{\infty}$ and $\|\nabla \mu\|_{\infty} \lesssim \| \gamma \|_{W^{1,\infty}}$.  
		
		We define $v_h:=e^{\mu/h}u'_h$, then on the boundary:
		$$h\partial_{n}u'_h+\gamma u'_h=h\partial_{n} v_h=0.$$
		Hence the function $v_h\in H^1(M',\R)$ is solution of:
		\begin{equation}
			\left\{\begin{array}{ll}
				-\frac{h^2e^{\mu/h}}{K'_h}\mathrm{div}_{g'}(K'_h \nabla_{g'})e^{-\mu/h}v_h +V_h' v_h =E_h v_h \mbox{ on } M',\\
				\partial_{n} v_h=0 \mbox{ on } \partial M'. 
			\end{array}\right.
		\end{equation}
		Note that:
		$$-h^2\frac{e^{\mu/h}}{K'_h}\mathrm{div}_{g'}(K'_h\nabla_{g'})e^{-\mu/h}v_h=-h^2\frac{e^{2\mu/h}}{K'_h}e^{-\mu/h}\mathrm{div}_g(e^{\mu/h}e^{-2\mu/h}K'_he^{\mu/h}\nabla_{g'})e^{-\mu/h}v_h$$
		$$=-h^2\frac{e^{2\mu/h}}{K'_h}\mathrm{div}_{g'}(K'_he^{-2\mu/h}\nabla_{g'})v_h+(\frac{h}{K'_h}\mathrm{div}_{g'}(K'_h\nabla_{g'})\mu-|\nabla_{g'}\mu|_{g'^{-1}}^{2})v_h.$$
		So if we put $\kappa_{h} :=e^{-2\mu /h}K'_{h}$ and $\mathcal{V}_{h} := h (K'_{h})^{-1}\mathrm{div}_{g'}(K'_{h}  \nabla_{g'})\mu-|\nabla_{g'}\mu|_{g'^{-1}}^{2}+V'_h$, then
		
		\begin{equation}
			\left\{\begin{array}{ll}
				\left(-\frac{h^2}{\kappa_h}\mathrm{div}_{g'}(\kappa_h\nabla_{g'}) +\mathcal{V}_h\right) v_h=E_h v_h  \mbox{ on } M',\\
				\partial_{n} v_h=0 \mbox{ on } \partial M'.
			\end{array}\right.
		\end{equation}
		
		We use again the proposition to get a covering of $\partial M' \subset \bigcup_{j=1}^N U'_j$ (seen as a subset of $\overline{M'}$) and $W^{2,\infty}$ diffeomorphisms $\Psi'_j: U'_j \to \R^2=\R_x\times \R_y$ 
		such that $$\Psi'_j(U'_j)\subset B(0,1)\cap \{y\geq 0\}$$ and the metric takes the
		form \begin{equation}		
			b'_{|\partial \Omega}(x,y)=\begin{pmatrix}
				b'(x) & 0\\
				0 & 1
			\end{pmatrix}
		\end{equation}
		on the boundary with $b'\in W^{1,\infty}$.
		
		We now perform the second gluying by defining a covering of $\partial M'$ (seen as a subset of $\tilde{M}$ ):
		$$\partial M'\subset \bigcup_{j=1}^{N}U'_j\times\{\pm i\}=:\bigcup_{j=1}^{N}\tilde{U}_j$$
		where we identify the points in $\left(U'_j\cap \partial M'\right)\times \{\pm i\}$, and define the map
		$$\tilde{\Psi}_j:z=(x,\epsilon)\in U'_j\times \{\pm i\}\to \left\{
		\begin{array}{ll}
			\Psi'_j(x) & \mbox{ if } \epsilon=1, \\
			S\circ \Psi'_j(x) & \mbox{ if } \epsilon=-1,
		\end{array}
		\right.$$
		
		For the same precedents reasons, $(\tilde{\Psi}_j,\tilde{U}_j)$ form a $W^{2,\infty}$ structure on $(\tilde{M}, \tilde{g})$ where $\tilde{g}$ Lipschitz such that   $\tilde{g}(x,\pm i)=g'(x)$. We define the extended weight $\tilde{K}_h$ by $\tilde{K}_h(x,\pm i):= \kappa_h(x).$ For same reason use previously,  $\tilde{K}_h$ is continuous on $\partial M'$ and Lipschitz outside, thus  $\tilde{K}_h\in W^{1,\infty}(\tilde{M},\R)$  by the jump formula. 
		
		We now extend $v_h$ on $M'$ by $\tilde{u}_h(x,\pm i)=v_h(x)$ and $\tilde{V}_h(x,\pm i)=\mathcal{V}_h(x)$.
		
		 Away from  $\partial M'$ it is clear that 
		 \begin{equation}
		 	\label{eq: EDP_weight_boundarycase}
		 	-\frac{h^2}{\tilde{K}_h}\mathrm{div}_{\tilde{g}}(\tilde{K}_h\nabla_{\tilde{g}})\tilde{u}_h +\tilde{V}_h \tilde{u}_h= E_h \tilde{u}_h. 
		 \end{equation}
		 Near a point $x\in \partial M'$ we can work in a coordinate chart $(\tilde{\Psi}_j, \tilde{U}_j)$. In this coordinate chart, the function
		$\tilde{u}_h$ is defined by
		$$ \tilde{u}_h=\left\{
		\begin{array}{ll}
			v_h(x,y) &\mbox{ if } y>0,\\
			v_h(x,-y) &\mbox{ if } y<0.
		\end{array}
		\right.$$
		In $\{\pm y>0\}$, $\tilde{u}_h$ satisfies the equation \ref{eq: EDP_weight_boundarycase}, and near $\partial M'$ in our coordinate systems, we have
		$$\tilde{u}_h(x,y)=v_h(x,y)\,1_{y>0}\,+v_h(x,-y)\, 1_{y<0},$$ and we have
		$$\partial_x \tilde{u}_h(x,y)=\partial_x (v_h(x,y))\,1_{y>0}\,-\partial_x(v_h(x,-y))\, 1_{y<0};$$
		$$\partial_y \tilde{u}_h(x,y)=\partial_y (v_h(x,y))\,1_{y>0}\,-\partial_y(v_h(x,-y))\, 1_{y<0}\,+(v_h(x,0^{+})-v_h(x,0^{-}))\otimes \delta_{y=0}$$
		$$=\partial_y (v_h(x,y))\,1_{y>0}\,-\partial_y(v_h(x,-y))\, 1_{y<0}.$$
		So the extension define a function in $H^1(\tilde{M},\R)$. Then with the form of $\tilde{g}$ on $\partial M'$ and by the jump formula:
		\begin{align}
			h^2\tilde{K}_h^{-1}\mathrm{div}_{\tilde{g}}(\tilde{K}_h\nabla_{\tilde{g}})\tilde{u}_h(x,y)&= (\tilde{V}_h(x,y)-E_h)\tilde{u}_h(x,y)\nn\\
			&+\left(\partial_y v_h(x,0^{+})+\partial_y v_h(x,0^{-})+0\times (\partial_x v_h(x,0^{+})-\partial_x v_h(x,0^{-})) \right)\otimes \delta_{y=0}\nn\\
			&=(\tilde{V}_h(x,y)-E_h)\tilde{u}_h(x,y)
		\end{align}
		because $\partial_y v_h=0$ on $\partial M'$.
		
		Finally, we have obtain a function $\tilde{u}_h\in H^1(\tilde{M},\R)$ with $\tilde{M}$ closed $W^{2,\infty}$ Riemannian surface with Lipschitz metric $\tilde{g}$. Moreover, $\tilde{u}_h$ is solution of:
		$$-\frac{h^2}{\tilde{K}_h}\mathrm{div}_{\tilde{g}}(\tilde{K}_h\nabla_{\tilde{g}})\tilde{u}_h +\tilde{V}_h \tilde{u}_h= E_h \tilde{u}_h$$
		 on $\tilde{M}$ with $\tilde{K}_h=e^{\tilde{\phi}_h/h}$ with $\tilde{\phi}_h=\tilde{\phi}+\tilde{\mu}+\mathcal{O}(h)\in W^{1,\infty}(\tilde{M},\R) $, $\tilde{V}_h=\tilde{V} +\mathcal{O}(h)\in L^\infty(\tilde{M},\R)$ and $E_h=E+\mathcal{O}(h)\in \R$.
		 Thus we can apply the Theorem \ref{thm: closed_Riem_surf} and conclude the proof.
	\end{proof}
	
	\appendix
	\section{Section Preliminaries proofs}
	
	We proof here the lemmas from Preliminaries.
	
	\begin{lem}
		\label{lem: rég_C1_proof}
		If $P_{h}\, u_{h}=0$ on an open set $\Omega$ of $\R^2$ in the sense of distributions and $u_{h}\in H_{loc}^1(\Omega)$, then $u_{h} \in \mathcal{C}^1(\Omega)$.
	\end{lem}
	
	\begin{proof}
		Firstly, as a function in $H_{loc}^1(\Omega)$,  solution of an elliptic equation of $$\left(-\frac{h^2}{K_h}\mathrm{div}(K_h\nabla)+ V_h\right)u_h=0$$
		with $K_h\in W^{1,\infty}$ and $V_h\in L^\infty$, $u_h$ is in  $ H_{loc}^2(\Omega)$ by elliptic regularity ( see Theorem $8.8$ of \cite{Gilbarg2001}). Thus it imply that $u_h$ is locally Lipschitz by Sobolev embedding. Therefore the function $f_h:=-h^2K_h^{-1}\nabla K_h\nabla u_h +V_h u_h$ define a $L_{loc}^{\infty}$ function.
		Let $\Omega'\subset \subset \Omega$ bounded open set, we define $\Gamma(x):= -\frac{1}{2\pi} \log |x|$  the fundamental solution of the $2D$ Laplacian on the plane.
		the function define by convolution $$\tilde{u}_h:= \Gamma \star f_h$$ is well defined on $\Omega'$ and verify that $$-h^2\Delta (u_h-\tilde{u}_h)=0$$ in the sense  of distributions. Consequently, $u_h-\tilde{u}_h$ is harmonic by the Weyl's lemma and so is smooth. So we just need to verify that $\tilde{u}_h\in \mathcal{C}^{1}(\Omega')$. For that we identified $\R^2$ to $\C$. With a first order expansion, we see that
		$$\Gamma(x+t)-\Gamma
		(x)=\frac{1}{2\pi}\log\left| 1+\frac{t}{x} \right|=\Re\frac{t}{2\pi x}+ \mathcal{O}\left(W\left(\frac{x}{t}\right)\right)$$
		where $W(\xi)= \begin{cases}
			\frac{1}{|\xi|^2 } & |\xi|>2; \\
			\frac{1}{|\xi|}+\log\left|1+\frac{1}{\xi}\right| & |\xi| \leq 2.
		\end{cases}$
		
		Thus we have the following expansion  $$\tilde{u}_h(x+t)=\tilde{u}_h(x)+\Re \left(t\,\frac{1}{2\pi \cdot }\star f_h(x)\right)+ \mathcal{O}(||f_h||_{L^\infty} ||W(\cdot/t)||_{L^1}).$$ The second term in the expansion is a linear application from $t\in\C$ to $\R$ which depend continuously of $x$ and $$ ||W(\cdot/t)||_{L^1(\Omega')}\leq  |\Omega'|^{1-1/q}||W(\cdot/t)||_{L^q(\Omega')}\leq |\Omega'|^{1-1/q}t^{2/q}||W||_{L^q(\R^2)}$$ where $1<q<2$. It is clear that $||W||_{L^q(\R^2)}<\infty$, thus 
		$$\tilde{u}_h(x+t)=\tilde{u}_h(x)+\Re \left(t\,\frac{1}{2\pi }\star f_h(x)\right)+ o(t)$$ with $\frac{1}{2\pi \cdot }\star f_h\in \mathcal{C}^0(\Omega')$. Hence the $\mathcal{C}^1$ regularity.
		
	\end{proof}
	
\begin{lem}
	Let $\Omega_h$ be an open set of $B(0,R_0)$ with Poincaré constant $k(h)<<h$. 
	There is a weak solution (in the sense of distributions) $u\in H_0^1(\Omega_h,\R)$ to the equation
	\begin{equation}
		h^2K_h^{-1}\mathrm{div}(K_h\nabla u)=v
	\end{equation}
	for all $v\in L^2\cap L^\infty(\Omega_h,\R )$: for all $\chi \in \mathcal{C}_0^\infty(\Omega_h,\R)$ $$\int_{\Omega_h}h^2K_h\nabla u\cdot \nabla\chi=\int_{\Omega_h} K_hv\chi.$$
	The solution $u$ verify	
	\begin{itemize}
		\item $||u||_{H^{1,k(h)}(K_hdx)}\leq A\frac{k(h)^2}{h^2}||v||_{L^2(K_hdx)}$,
		\item $||u\, \sqrt{K_h}||_{L^\infty}\leq  A\frac{k(h)^2}{h^2}||v\, \sqrt{K_h}||_{L^\infty}$.
	\end{itemize}
	with $A>0$ an absolute constant and where $$||\cdot||_{H^{1,k(h)}(K_hdx)}:= \sqrt{\int \left(|\cdot|^2+|k(h)\nabla \cdot|^2 \right)K_hdx}.$$ 
\end{lem}
	
	\begin{proof}
		for the construction of $u$, we search it has a minimum of the following functional: 
		\begin{equation}
			\Phi(w)=\int \frac{h^2K_h}{2}|\nabla w|^2 +K_hvw
		\end{equation}
		with $w\in H_0^{1}(\Omega_h)$.
		For all $w\in H_0^{1}(\Omega_h)$:
		\begin{equation}
			\Phi(w)\geq\frac{ h^2}{2} ||\sqrt{K_h}\nabla w||_{L^2}^2-k(h) ||\sqrt{K_h}v||_{L^2}||\sqrt{K_h}\nabla w||_{L^2}\\
			\geq -\frac{k(h)^2}{2h^2}||\sqrt{K_h}v||_{L^2}^2.
		\end{equation}
		Thus this functional is bounded from below. Note that $\Phi(u)>0=\Phi(0)$ if $||k(h)\sqrt{K_h}\nabla u||_{L^2}>\frac{2k(h)^2}{h^2}||\sqrt{K_h}v||_{L^2}$. Then for any minimizing sequence $u_k\in \mathcal{C}_0^\infty$:
		
		\begin{equation}
			\frac{\Phi(u_{k})+\Phi(u_{k'})}{2}-\Phi\left(\frac{u_{k}+u_{k'}}{2}\right)=\frac{1}{8}h^2||\sqrt{K_h}\nabla(u_k-u_{k'})||_{L^2}.
		\end{equation} 
		Therefore $u_k$ is a Cauchy sequence which will converge to $u\in H_0^{1}(\Omega_h)$.
		As a minimizer, $u$ will be  solution of 
		\begin{equation}
			\int h^2 K_h\nabla u \nabla \phi +K_hv\phi=0
		\end{equation}
		for all $\phi \in \mathcal{C}_0^\infty(\Omega_h)$. So $u$ is a weak solution of $h^2K_h^{-1}\mathrm{div}(K_h\nabla u)=v$ and 
		\begin{equation}
			||k(h)\sqrt{K_h}\nabla u||_{L^2}\leq 2\frac{k(h)^2}{h^2}||\sqrt{K_h}v||_{L^2}.
		\end{equation}
		Then by the Poincaré inequality:
		\begin{equation}
			||\sqrt{K_h}u||_{L^2}\leq k(h) ||\sqrt{K_h}\nabla u||_{L^2} \leq  2\frac{k(h)^2}{h^2}||\sqrt{K_h}v||_{L^2}.
		\end{equation}
		
		To show the inequality on $L^\infty-$norm, we use a Di Giorgi method. The proof consist of three steps. For convenience, we change the definition of the functional $\Phi$ by $$ \Phi(u)=\int K_h|\nabla u|^2 +K_hv_h u$$ with $v_h=\frac{2v}{h^2}$.
		
		\begin{lem}($\mbox{Step } 1$)\\
			Assume that the Poincaré constant $k(h)$ is small. Then any smooth $\delta-$minimizer $u_{\delta}$ of $\Phi$ (ie for $u\in C_0^\infty(\Omega_h)$, $\Phi(u)\geq \Phi(u_{\delta})-\delta$) satisfies
			\begin{equation}
				\int_{B\cap \Omega_h} u_{\delta}^{2} \leq Ck(h)^2(k(h)^2||v_h\sqrt{K_h}||_{\infty}+\delta)
			\end{equation}
			for every unit ball $B \subset \mathbb{R}^{2}$.
		\end{lem}

		\begin{proof}
			This proof relies on Poincaré inequality.
			
			Let $\varphi$ be a smooth positive radial function such that 
			\begin{itemize}
				\item $\varphi=1$ in $B(0,1)$, 
				\item $0<\varphi \leq 1$,
				\item $|\nabla \varphi| \leq  \varphi$.
			\end{itemize}
			
			Applying the Poincaré inequality to $\varphi u\sqrt{K_h}$, we get
			\begin{align}
				&\int \varphi^{2} u_{\delta}^{2}K_h\leq k^{2} \int|\varphi \nabla u_{\delta}+u_{\delta} \nabla \varphi+u_{\delta}\varphi \nabla \phi_h/h|^{2}K_h\nn\\
				&\leq 3 k(h)^{2}\left(\int \varphi^{2}|\nabla u_{\delta}|^{2}K_h+\int u_{\delta}^{2}|\nabla \varphi|^{2}K_h+\int \varphi^{2} u_{\delta}^2|\nabla \phi_h/h|^{2}K_h\right)\nn\\ 
				&\leq 3k(h)^{2} \int \varphi^{2}|\nabla u_{\delta}|^{2}K_h+3 k(h)^{2}(1+||\nabla \phi_h/h||_{\infty}^2) \int u_{\delta}^{2} \varphi^{2}K_h\nn\\
			\end{align}
			
			Thus 
			\begin{equation}
				\label{eq: toollemlocsoldivtypeq}
				\int \varphi^{2} u_{\delta}^{2}K_h\leq \frac{3k(h)^2}{1-C(k(h)/h)^2}\int \varphi^{2}|\nabla u_{\delta}|^{2}K_h\leq 4k(h)^2\int \varphi^{2}|\nabla u_{\delta}|^{2}K_h
			\end{equation}
			for $k(h)$ small enough.

			Now we consider the function $\tilde{u}_{\delta}=(1-\varphi^2) u_{\delta}$.
			We have
			\begin{align}
				&\Phi(\tilde{u}_{\delta})=\int K_h|(1-\varphi^2) \nabla u_{\delta}-u_{\delta} \nabla \varphi^2|^{2}+\int K_hv_h(1-\varphi) u_{\delta}\nn\\
				& \leq \int K_h(1-\varphi^2)^{2}|\nabla u_{\delta}|^{2}+ \int 2K_h|\nabla u_{\delta}| | u_{\delta}||\nabla \varphi^2|+\int K_h u_{\delta}^{2}|\nabla \varphi|^{2}+K_h\int v_h(1-\varphi^2) u_\delta\nn\\
				&\leq \int K_h|\nabla u_{\delta}|^{2}+\int K_hv_h u_{\delta}-\int K_h|\nabla u_{\delta}|^{2}\varphi^2+4 \int  K_h|\nabla u_{\delta}| | u_{\delta}| \varphi^2 +4 \int K_h u_{\delta}^{2} \varphi^2-\int K_hv_h \varphi^2 u_{\delta}\nn.
			\end{align}
			
			(we used the fact, $|\nabla \varphi^2| \leq 2 \varphi^2$ and $ \varphi \leq 1$).
			
			Since $u_{\delta}$ is a $\delta-$minimiser, we must have
			
			\begin{equation}
				\int K_h|\nabla u_{\delta}|^{2} \varphi^2 \leq 4\left(\int K_h|\nabla u_{\delta}| |u_{\delta}| \varphi^2+\int K_h u_{\delta}^{2} \varphi^2\right)+\delta+\int K_h| v_h \varphi^2 u_{\delta} |.
			\end{equation}
			
			By the equation \ref{eq: toollemlocsoldivtypeq}, we know that
			\begin{equation}
				\int K_h\varphi^{2} u_{\delta}^{2}\leq 4k(h)^2\int K_h\varphi^{2}|\nabla u_{\delta}|^{2}
			\end{equation}
			and
			\begin{equation}
				\int K_h|\nabla u_{\delta}| |u_\delta| \varphi^2 \leq \sqrt{\int K_h u_{\delta}^{2} \varphi^2} \sqrt{\int K_h | \nabla u_{\delta}|^{2} \varphi^2} \leq 2 k(h) \int K_h|\nabla u_\delta|^{2} \varphi^2.
			\end{equation}
			
			So for sufficiently small $k(h)$,
			$$ 4\left(\int K_h|\nabla u_{\delta}| |u_\delta| \varphi^2+\int K_hu_{\delta}^{2} \varphi^2\right) \leq \frac{1}{2} \int K_h|\nabla u_{\delta}|^{2} \varphi^2.$$
			Hence
			\begin{multline}
				\int K_h|\nabla u_{\delta}|^{2} \varphi^2\leq 2\delta+ 2\int K_h|v_h u_{\delta} \varphi^2| \\
				\leq 2\left(\delta+||v_h\sqrt{K_h}||_{\infty}\sqrt{\int \varphi^2} \sqrt{\int u_{\delta}^{2} \varphi^2K_h} \right) \leq C\left(\delta+ k(h) ||v_h\sqrt{K_h}||_{\infty} \sqrt{\int K_h|\nabla u_{\delta}|^{2} \varphi^2}\right).
			\end{multline}
			
			Thus with the equation \ref{eq: toollemlocsoldivtypeq}, it follows  that 
			\begin{equation}
				\int u^{2} \varphi^2\leq Ck(h)^2\left(\delta +k(h)^2||v_h\sqrt{K_h}||_{\infty}^2\right).
			\end{equation}
		\end{proof}

		\begin{lem}($\mbox{Step } 2$)
			
			Assume $\Omega$ is a set such that
			$$m_2(\Omega \cap Q)\leq k^2 \ll 1$$ 
			for any unit square $Q$, then 
			the Poincaré constant of $\Omega$ is smaller than $Ck$. 
		\end{lem}
		
	\begin{proof}
		It is clear that if we can show that $$\int_{Q} |f|^2\leq Ck^2 \int_{Q}|\nabla f|^2$$
		for  any $f\in \mathcal{C}_0^\infty(\Omega)$ and unit cube $Q$, then we will obtain the inequality  by tiling the plane with unit squares.
		Let $f\in \mathcal{C}_0^\infty( \Omega)$, $x\in \Omega \cap Q$, we have:
		$$f(x)=\int_{[x_0,x]} \nabla f$$ 
		for all $x_0\in Q\setminus \Omega$. Hence
		$$|f(x)|\leq \frac{1}{m_2(Q\setminus \Omega)}\int_{Q\setminus \Omega}\int_{[x_0,x]}|\nabla f|dx_0   $$ Let note $B(x,\sqrt{2})$ the ball centre on $x$ with radius $\sqrt{2}$. Then $$|f(x)|\leq \frac{1}{m_2(Q\setminus \Omega)}\int_{B(x,\sqrt{2})}\int_{[y(x_0),x]}|\nabla f| 1_Qdx_0$$ where $y(x_0)$ is the projection of $x_0$ on the boundary of $B(x,\sqrt{2})$ by the radius passing through $x_0$. Thus $$|f(x)|\leq \frac{1}{m_2(Q\setminus \Omega)} \int_{0}^{\sqrt{2}}\int_{0}^{\sqrt{2}}\int_{0}^{2\pi} |\nabla f(t,\theta)|1_Q(r,t,\theta) d\theta \,dt\,  rdr$$ $$\leq \frac{\sqrt{2}}{m_2(Q\setminus \Omega)} \int_{0}^{\sqrt{2}}\int_{0}^{2\pi} |\nabla f(t,\theta)|1_Q(t,\theta)d\theta\, dt $$ $$\leq \left(\int_{0}^{\sqrt{2}}\int_{0}^{2\pi} |\nabla f(t,\theta)|^{3/2}1_Q(t,\theta) d\theta \,tdt  \right)^{2/3}\left(\int_{0}^{\sqrt{2}}\int_{0}^{2\pi} t^{-1/2}d\theta\, dt \right)^{1/3}$$ $$\leq \frac{C}{m_2(Q\setminus \Omega)}\left(\int_Q |\nabla f(x_0)|^{3/2}dx_0 \right)^{2/3}\leq \frac{C'}{m_2(Q\setminus \Omega)}\left(\int_Q |\nabla f(x_0)|^{2}dx_0 \right)^{1/2}.$$ Hence if $k^2\leq 1/2$ then:
		\begin{equation}
			\int_Q |f(x)|^2dx=\int_{Q\cap \Omega} |f(x)|^2dx\leq \tilde{C}k^2\int_Q|\nabla f(x_0)|^{2}dx_0.
		\end{equation}
		
	\end{proof}
		
		$\mbox{Step }3$:
		As done by \cite{logunov2020landisconjectureexponentialdecay}, we will now run a Di Georgi scheme. For that, we need a special smooth function:
		Let $w\in C_0^\infty(\Omega_h)$. Fix $t>0$ a level  and $\mu>0 $ an extremely small number.
		Let $\Theta$ be a $C^\infty$-smooth function on $\mathbb{R}$ such that:
		
		\begin{itemize}
			\item $ 0 \leq  \Theta' \leq 1,$
			\item $ \Theta(0)= 0,$
			\item $ x-4\mu\leq \Theta(x) \leq x,$
			\item $\Theta(x) = t-2\mu \text{ on } (t -\mu, t+\mu).$
			
		\end{itemize}
		Let $\tilde w = \Theta \circ w$. Then $$|w-\tilde w| \leq 4\mu \mbox{ and } 
		|\nabla \tilde w| \leq |\nabla w|$$ pointwise.
		Thus, if $w$ was an $\delta$-minimizer,
		then $\tilde w$ is an $\delta+ A_h||v_h\sqrt{K_h}||_\infty\mu$-minimizer with $A_h=4||\sqrt{K_h}||_\infty|\Omega_h|$ and $\mu$ can be choose arbitrary small.
		Define $$\Theta_{-}(x)=
		\begin{cases}
			\Theta(x), \quad x\leq t+\mu\\
			t-2\delta, \quad x\geq t+\mu
		\end{cases}
		\quad \text{and} \quad \Theta_+=\Theta-\Theta_-.  $$
		The function $\tilde w$ splits into two smooth compactly supported functions:
		$$ \tilde w= \tilde w_{-}+ \tilde w_{+}, \text{ where } \tilde w_\pm = \Theta_\pm \circ w. $$
		The function $\tilde w_{+}$ is compactly supported in $\{w>t\}$ ($\text{supp } \tilde w_+ \subset\{w \geq t+\mu\}$)
		and $\nabla \tilde w_{+}$ and $\nabla \tilde w_{-}$ have disjoint supports.
		Hence: 
		$$ \Phi(\tilde w)=  \Phi(\tilde{w}_-)+ \int_{\{w>t\}}|\nabla \tilde{w}_+ |^2  + \int_{\{w>t\}} v_h \tilde{w}_+,$$
		Thus if we replace $ \tilde w_+$ by any smooth function supported on $\{w>t\}$, we can see that $\tilde{w}_+$ is an $(\delta+A_h||v_h\sqrt{K_h}||_\infty\mu)$-minimizer of the functional $\Phi$ on the new domain $\{w> t \}$.
		
		Now, if $u_{\delta}$ is an $\delta$-minimizer  with $\delta= k(h)^2||v_h\sqrt{K_h}||_{\infty}^2\delta_0$, then
		$$ \int_{B\cap \Omega_h} u_{\delta}^2K_h \leq Ck(h)^4||v_h\sqrt{K_h}||_{\infty}^2(1+\delta_0)$$ for every unit ball $B$ by the step $1$. So 
		$$m_2(\{u_{\delta}\sqrt{K_h}>t_0\} \cap B)\leq C\frac{k(h)^4||v_h\sqrt{K_h}||_{\infty}^2(1+\delta_0)}{t_0^2}.$$

		Choose $t_{0}=C' k(h)||v_h\sqrt{K_h}||_{\infty}$ with sufficiently large constant $C'$. Then the domain
		$\Omega_{1,h}=\left\{u_{\delta}>t_{0}\right\}$ satisfies $$m_2(\Omega_{1,h}\cap B) \leq \frac{Ck(h)^2}{C'} \left(1+\delta_0\right)$$ for any unit ball $B$. By the Step $2$,
		the Poincaré constant of $\Omega_{1,h}$ is at most $$k(h)_1=k(h)k_1:=k(h)\left(\frac{1+\delta_0}{4}\right)^{1/2}.$$
		Also $u_{{\delta_1}}=\tilde{u_{\delta}}_{+} \in C_{0}^{\infty}\left(\Omega_{1,h}\right)$ will be an $k(h)^2||v_h\sqrt{K_h}||_\infty^2\delta_{1}$-minimizer of $\Phi$ 
		in $C_{0}^{\infty}\left(\Omega_{1,h}\right)$ with $\delta_1=(\delta_0+A_h\mu_0)$ where $\mu_0>0$ can be chosen arbitrary small. Finally, note
		that $$u_{\delta}\sqrt{K_h} \leq t_{0}+u_{{\delta_1}}\sqrt{K_h}+4k(h)^2||v_h\sqrt{K_h}||_\infty\mu_{0}$$ everywhere in $\Omega_h$.
		
		For coherence of notations, we will note $\Omega_{0,h}=\Omega_h$, $u_{\delta_0}=u_{\delta}$ and $k_0=1$.
		We can now repeat this construction with $ u_{\delta_1}$, $\Omega_{1,h}$, $\delta_1$ instead of $u_{\delta_0}$, $\Omega_{0,h}$, $\delta_0$ to
		get $u_{\delta_2}$, $\Omega_{2,h}$, $\delta_{2}$ and so on. We shall get a sequence of domains $\Omega_{j,h}$, functions
		$u_{\delta_j} \in C_0^{\infty}(\Omega_{j,h})$ and numbers $k_j,t_j,\delta_j, \mu_j>0$ such that 
		
		$$m_2(\Omega_{j,h}\cap B) \leq Ck(h)^4||v_h\sqrt{K_h}||_\infty^2 k_j^2(k_j^2 +\delta_{j}) \quad \text{ for any unit ball } B,  $$
		$$ t_j=C'k(h)\sqrt{k_j}||v_h\sqrt{K_h}||_{\infty},  \quad  \Omega_{j+1,h} \subset \Omega_{j,h},$$
		$$ \delta_{j+1}= \delta_{j}+ A_h\mu_{j}, \quad k_{j+1}^2=k_j\frac{(k_j^2+\delta_{j})}{4},  $$
		
		and 
		$$ u_{\delta_0}\sqrt{K_h} \leq t_0+t_1+...+t_{j}+u_{j+1}\sqrt{K_h}+4k(h)^2||v_h\sqrt{K_h}||_\infty(\mu_0+\mu_1+...+\mu_{j}) \quad \text{for any } j\geq 0.$$
		
		In this construction, we can choose $\mu_j>0$ as  small as we want, so putting $\mu_j=\frac{c_h\delta_0}{2^j }$ with sufficiently small absolute constant $c>0$, we can guarantee that all $\delta_j \leq 2 \delta_0$.
		
		If $k_j^2 \geq 2\delta_0$, then $$k_{j+1}^2\leq\frac{k_j^3}{2}$$ so the sequence is strictly decreasing.
		
		Let $l$ be the first index for which $k_l^2<2\delta_0$. Then, we also have (for $C/C'$ small enough) 
		$m_2(\Omega_{l+1,h} \cap B)<\frac{k(h)^2\delta_0}{2}$ for every unit ball $B$, and for all $j\leq l$, we have
		$$k_j^2 \leq 2^{-j}.$$
		
		Whence $t_j \leq C'k(h)||v_h\sqrt{K_h}||_\infty \sqrt{2}^{-j}$ and 
		\begin{multline}
			u_{\delta}\sqrt{K_h} \leq t_0+t_1+...+t_{l}+4k(h)^2||v_h\sqrt{K_h}||_\infty(\mu_0+...+\mu_{l}) \\ \leq (C_1+C_2k(h))k(h)||v_h\sqrt{K_h}||_\infty=: C_0k(h)||v_h\sqrt{K_h}||_\infty
		\end{multline}
		in $\Omega_h\setminus \Omega_{l+1,h}$ because $u_{\delta_{l+1}}=0$ in $\Omega_h\setminus \Omega_{l+1,h}$.
		Thus
		$$m_2(\{ u_{\delta}\sqrt{K_h}>C_0k(h)||v_h\sqrt{K_h}||_\infty\}\cap B)<k(h)^2\delta_0 /2.$$
		Considering $-u_{\delta}$ instead of $u_{\delta}$, we conclude that also $$m_2(\{ u_{\delta}\sqrt{K_h}<-C_0k(h)||v_h\sqrt{K_h}||_\infty\}\cap B)<k(h)^2\delta_0/2$$ and
		therefore
		$$m_2(\{ |u_{\delta}\sqrt{K_h}|>C_0k(h)||v_h\sqrt{K_h}||_\infty\}\cap B)<k(h)^2\delta_0$$
		for every $k(h)^2||v_h\sqrt{K_h}||_\infty\delta_0$-minimizer $u_\delta$. Since the true minimizer $u$ is the limit
		of $k(h)^2||v_h\sqrt{K_h}||_\infty\delta_0$-minimizers in $L^2(\Omega_h)$, we get 
		$$ m_2(\{|u\sqrt{K_h}|>C_0k(h)||v_h\sqrt{K_h}||_\infty\}\cap B)=0 \quad \text{ for any unit ball } B,$$
		so $m_2(\{|u\sqrt{K_h}|>C_0k(h)||v_h\sqrt{K_h}||_\infty\})=0$. 
		Thus 
		\begin{equation}
			||u\sqrt{K_h}||_\infty\leq C_0k(h)||v_h\sqrt{K_h}||_\infty \mbox{ a.e}.
		\end{equation}	
		In order to win one more $k(h)$, we will use a scaling: Let be $\tilde{\Omega}_h=\frac{k_0}{k(h)} \Omega_h$, $\tilde{u}=\frac{k_0^2}{k(h)^2}u(\frac{k(h)}{k_0}\cdot )$ and $\tilde{v}_h= v_h(\frac{k(h)}{k_0}\cdot )$. So $$\tilde{\Phi}(\tilde{u})=\frac{k_0^4}{k(h)^4}\Phi(u),$$ hence $\tilde{u}$ is a minimizer of $\tilde{\Phi}$ and if $k_0$ is small enough, we can use again the above result to obtain:
		
		\begin{equation}
			||u\sqrt{K_h}||_\infty\leq \tilde{C}_0k(h)^2||v_h\sqrt{K_h}||_\infty= 2\tilde{C}_0\frac{k(h)^2}{h^2}||v\sqrt{K_h}||_\infty \mbox{ a.e}.
		\end{equation} 
	\end{proof}

	\begin{prop}
		Let $u_h$ be a $\mathcal{C}^1$ solution of $P_{h}u_h=0$ on a ball $B(x,r_h)$ of $B(0,R_0)$ where $r_h<r_{0,h}=r_0h$ with $r_{0}$ a small constant independent of $h$ and $u_h$ for $h$ small enough. If $u_h>0$ on $\partial B(x,r_h)$, then $u_h>0$ on $B(x,r_h)$. 
	\end{prop}
	
	\begin{proof}
		To prove this result, we will need a auxiliary lemma: 
		
			\begin{lem}
			\label{lem: C01_imply_H01}
			Let $\Omega_h$ be a bounded open set and let a function $f\in C^1(\overline{\Omega_h})$ be zero on $\partial \Omega_h$. Then $f \in H_{0}^{1}(\Omega_h)$.
		\end{lem}
		
		\begin{proof}
			This proof is taken from lemma $6.13.$ of \cite{logunov2020landisconjectureexponentialdecay}.
			
			Denote by $\Omega_{\epsilon,h}$ the set of points $x$ in $\Omega_h$ with distance to the boundary of $\Omega_h$ at least $\epsilon>0$. Let $\eta$ be a function in $\mathcal{C}_{0}^{\infty}(\Omega_h)$ with the following properties: 
			\begin{itemize}
				\item $\eta=1$ on $\Omega_{\epsilon,h}$,
				\item $0\leq \eta \leq 1$ and $|\nabla \eta|\leq C\epsilon$ in $\Omega_h$.
			\end{itemize}
			The function $f\eta$ is in $\mathcal{C}_0^1(\Omega_h)\subset H_0^1 (\Omega_h)$. We want to show that $f\eta$ converge  to $f$ in $H_0^1(\Omega_h)$. Observe that $|\nabla f|$ is uniformly bounded in  $\Omega_h$ by some constant $A(f)$, so $|f(x)|\leq A(f)\epsilon$ if the distance from $x$ to $\partial \Omega_h$ is smaller than $\epsilon$ and therefore 
			$$ |\nabla(f\eta)|\leq |\nabla f | |\eta| + |f||\nabla \eta| \leq A(f) + A(f)C \mbox{ in } \Omega_h \setminus \Omega_{\epsilon,h}.$$  Then 
			$$ \int_{\Omega_h}  |f - f\eta|^2 =  \int_{\Omega_h\setminus \Omega_{\epsilon,h}}  |f - f\eta|^2 \leq \int_{\Omega_h\setminus \Omega_{\epsilon,h}}  |f|^2 \to 0$$  and  $$\int_{\Omega_h}  |\nabla f - \nabla(f\eta)|^2 =  \int_{\Omega_h\setminus \Omega_{\epsilon,h}} |\nabla f - \nabla(f\eta)|^2 \leq 2\int_{\Omega_h\setminus \Omega_{\epsilon,h}}(|\nabla f |^2 + |\nabla(f\eta)|^2) \leq m_2(\Omega_h\setminus \Omega_{\epsilon,h})C_1A(f)^2.$$  Since $m_2(\Omega_h\setminus \Omega_{\epsilon,h}) \to 0$ as $\epsilon \to 0$, we conclude that $f\in H_0^1(\Omega_h)$.
		\end{proof}

		We can assume that $u_h$ is greater than a positive constant $\delta_h$ on $\partial B(x, r_h)$. Consider the set $\Omega_h = \{x \in B(x, r_h) : u_h(x) <\delta_h/2\}$. This is an open set strictly inside $B(x, r_h)$ and if $u_h$ is not positive in $B(x, r_h)$, then $\Omega_h$ is not empty. Since $u_h \in \mathcal{C}^1(\Omega_h)$ according to the lemma \ref{lem: rég_C1_proof} and $u_h = \delta_h/2$ on $\partial \Omega_h$, we know by lemma \ref{lem: C01_imply_H01} that $(u_h-\delta_h/2)\in H_0^{1}(\Omega_h)$.
		Since $\Omega_h \subset B(x,r_h)$, $\Omega_h$ has Poincaré constant smaller than $Cr_h^2$. 
		For $r_h$ small enough, according to the proposition \ref{prop: constrsmallsolSchrod}, there exists $\varphi_h=\tilde{K}_h^{-1/2}+\Tilde{\varphi}_h$ such that $(P_{h}-E_h)\,\varphi_h=0$ on $\Omega_h$ with 
		$\Tilde{\varphi}_h \in H_0^{1}(\Omega_h)$, $\tilde{K_h}\in \mathcal{C}^\infty(B(0,R_0),\R_+^{*})$ and $\frac{1}{2\sqrt{K_h}}\leq \varphi_h$. 
		Then the function $g_h=(\delta_h\tilde{K}_h^{1/2}\varphi_h/2 -u_h)$ is in $H_0^{1}(\Omega_h)$ and is a solution of $ (P_{h}-E_h) g_h=0$. 
		For any function $\beta_h \in \mathcal{C}_0^\infty(\Omega_h)$, 
		\begin{equation}
			h^2\int_{\Omega_h}K_h\nabla g_h \nabla \beta_h =-\int_{\Omega_h} K_h(V_h-E_h) g_h \beta_h
		\end{equation}
		and taking the limit in $H_0^{1}(\Omega_h)$, we get:
		\begin{multline}
			h^2\int_{\Omega_h}K_h|\nabla g_h|^2  \leq \Vert (V_h-E_h) \Vert_{L^\infty} \int_{\Omega_h} K_h g_h^2\\
			\leq \Vert (V_h-E_h) \Vert_{L^\infty}Cr_h^2\int_{\Omega_h} |\nabla \sqrt{K_h}g_h|^2 \leq 2\Vert (V_h-E_h) \Vert_{L^\infty}Cr_h^2\left(\int_{\Omega_h} K_h |\nabla g_h|^2+ g_h^2 |\nabla \sqrt{K_h}|^2\right).
		\end{multline}
		By definition  $ |\nabla \sqrt{K_h}|\leq ||\nabla \phi_h||_{\infty}/2h$, thus for $r_{0}$ small enough:
		\begin{equation}
		h^2\int_{\Omega_h}K_h|\nabla g_h|^2\leq  4\Vert (V_h-E_h) \Vert_{L^\infty}Cr_h^2\int_{\Omega_h} K_h |\nabla g_h|^2.
		\end{equation} 
		for $r_h<h (4C\Vert (V_h-E_h) \Vert_{L^\infty})^{-1/2}$ , this last inequality is true only if $g_h=0$. So $u_h=\delta_h\tilde{K}_h^{1/2} \varphi_h/2\geq \delta_h\sqrt{\tilde{K}_h}/4\sqrt{K_h}>0$ on $\Omega_h$ and so on $B(x,r_h)$.
	\end{proof}
	
	 \section{Weak Poincaré lemma}
	
	\begin{lem}
		\label{lem: Poincarelem_proof}
		Let $F\in L_{loc}^p(B(0,1),\R^2)$ ($1\leq p<\infty$) a vector field such that $\mathrm{div}(F)=0$ in the sense of distribution. Then there is a function $u\in W_{loc}^{1,p}(B(0,1))$ such that $$F=(F_1,F_2):=(\partial_{x_2}u,-\partial_{x_1}u).$$
	\end{lem}
	
	\begin{proof}
		This proof is taken from section $6.5$ of  \cite{logunov2020landisconjectureexponentialdecay}.
		Let  $$K_\epsilon(z) = \frac{1}{\epsilon^2} K(|z|/\epsilon)$$  be a smooth approximation to identity with $\mathrm{Supp}(K_{\epsilon}) \subset B(0, \epsilon)$. Define  $$F_{\epsilon} = F \star K_{\epsilon} = (F_1 \star K_{\epsilon}, F_2 \star K_\epsilon)$$  in the smaller ball $B(0, 1 - \epsilon)$. Then $F_\epsilon$ is divergence free in $B(0, 1 - \epsilon)$. Indeed, if $f \in \mathcal{C}_0^\infty(B(0, 1 - \epsilon))$, then by Fubini’s Theorem 
		$$ \int  (F \star K_{\epsilon})\nabla f =  \int  F (K_\epsilon \star \nabla f ) =  \int  F \nabla (K_\epsilon \star f ) = 0.$$  So there is a smooth function $u_\epsilon$ such that $F_{\epsilon} = \nabla \times u_{\epsilon}$ in $B(0, 1 - \epsilon)$. Fix $\delta \in (0, 1)$. By the Lebesgue theory $F \star K_{\epsilon}$ converge to $F$ in $L^p(B(0, 1 - \delta))$. Thus $\nabla u_\epsilon$ is a Cauchy sequence in $L^p(B(0, 1 - \delta))$. Let add a constant to $u_\epsilon$ so that $\int_B(0,1-\delta) u_\epsilon = 0$. By the Poincaré–Wirtinger inequality $u_{\epsilon}$ is a Cauchy sequence in $L^p(B(0, 1 - \delta))$. Thus we can find a function $\tilde{u}_{\delta}$ such that $u_{\epsilon}$ converge to $\tilde{u}_{\delta}$ in $W^{1,p}(B(0, 1 - \delta))$ and $\nabla \times \tilde{u}_{\delta} = \lim_{\epsilon \to 0}\nabla u_{\epsilon} = (F_1, F_2)$. For any $\delta_1, \delta_2 \in (0, 1)$, the gradients of $\tilde{u}_{\delta_1}$ and $\tilde{u}_{\delta_2}$ are the same in $B(0, 1 - \max(\delta_1, \delta_2))$ and therefore $\tilde{u}_{\delta_1}-\tilde{u}_{\delta_2}$ is constant almost everywhere in $B(0, 1 - \max(\delta_1, \delta_2))$. Finally, let us modify $\tilde{u}_{\delta}$ by subtracting a constant so that $\int_{B(0,1/2)} \tilde{u}_{\delta} = 0$ for all $\delta < 1/2$. Then  $u$ is well-defined by $u = \tilde{u}_{\delta}$ in $B(0, 1 - \delta)$.
	\end{proof}
	
	\section{Theorem $5.3$ of \cite{logunov2020landisconjectureexponentialdecay}}
	\label{sect: thm_Toy_proof}
	
	\begin{thm} (\cite{logunov2020landisconjectureexponentialdecay}, Theorem $5.3$.)
		\label{thm: Toy_proof}
		Let $D_j$ be a collection of $100-$separated unit disks on $\mathbb{R}^2=\mathbb{C}$ such that 
		$0 \notin \cup 3D_j$. Let $R> 10^4$, $0<r\leq R/4$. 
		Consider any harmonic function $u$ in $B(0,R) \setminus \cup D_j$ such that $u$ does not change sign in
		$(5D_j\setminus D_j) \cap B(0,R)$ for every $j$. Assume  that 
		$$\sup \limits _{B(0,R-R/64) \setminus \cup 3D_j}|u| \geq e^{- N} \sup \limits_{B(0,R) \setminus \cup 3D_j}|u|.$$
		Then 
		\begin{equation}
			\sup \limits _{B(0,r) \setminus \cup 3D_j}|u| \geq  \left(\frac r R \right)^{ L (R + N) } \sup \limits_{B(0,R) \setminus \cup 3D_j}|u|,
		\end{equation}
		with some absolute constant $L>0$.
	\end{thm}
	
	\begin{proof}
		Without loss of generality, we can take $\sup \limits _{B(0,R-R/64) \setminus 
			\cup 3D_j}|u|=1$. Fix $k=[C(N+R)]$ with sufficiently large $C>0$ and assume that
			\begin{equation}
				\label{eq: wrong_hyp}
				\sup \limits _{B(0,r) \setminus 
					\cup 3D_j}|u| \leq \left(\frac {r}{R} \right)^{3k}.
			\end{equation}
		
		Consider the domain $$\Omega:=\{r/2<|z|< R - 1 \} \setminus  \cup (3D_j).$$ 
		Let $W_1$ be the connected component of $\partial \Omega$ that intersects  the circle $C(0,r/2)$.
		Note that each point of $W_1$ is either on  $C(0,r/2)$ or lies on some $\partial 3D_j$ that intersects $C(0,r/2)$.\\
		
		\textbf{Estimate on $W_1$.}
		By Harnack inequality and Cauchy estimate on $5D_j \subset B(0,R)$, we have
		\begin{enumerate}
			\item  $\max\limits_{\partial3 D_j }|u| \leq A \min\limits_{\partial 3D_j}|u|$,
			\item  $\max\limits_{\partial 3D_j}|\nabla u| \leq A \min\limits_{\partial 3D_j}|u|$
		\end{enumerate}
		with $A>0$ because $u$ doesn't change its sign on $5D_j $.
		Hence on $W_1 \setminus C(0,r/2) $, we have
		$$ |u|,|\nabla u| \leq A \sup_{B(0,r)\setminus \cup (3D_j)} |u| \leq A  \left(\frac {r}{R} \right)^{3k}. $$
		If $x \in C(0,r/2) \setminus \cup (3D_j)$, then either $x \in 4D_j$ for some $j$ or $B(x, \min(1,r/2)) \subset B(0,r) \setminus \cup (3D_j) $. In the first case $u$ does not change sign in $B(x,1)$ and 
		$$|\nabla u(x)| \leq A|u(x)| \leq A\sup_{B(0,r)\setminus \cup (3D_j)} |u|\leq A\left(\frac {r}{R} \right)^{3k}.$$
		In the second case, we have 
		$$ |\nabla u(x)| \leq\frac{A}{\min(1,r/2)}\sup_{B(0,r)\setminus \cup (3D_j)}|u| \leq \frac{A}{\min(1,r/2)} \left(\frac{r}{R} \right)^{3k}. $$
		Thus in all cases, if  $C$ in the definition of $k$ is large enough, we have
		$$ \max_{W_1}|u|, \max_{W_1}|\nabla u| \leq \frac{A}{\min(1,r/2)} \left(\frac {r}{R} \right)^{3k} \leq \left(\frac r R \right)^{2k} \quad $$
		because  $$\left( \frac{R}{r} \right)^{k} \geq 4^k > A$$
		and
		$$\left( \frac{R}{r} \right)^{k} \geq 4^{k-1}  \frac{R}{r}  \geq \frac{2A}{r}.$$
		
		Let $W_2$ be the connected component of $\partial \Omega$ that intersects the circle $C(0,R-1)$.
		Note that each point of $W_2$ is either on  $C(0,R-1)$ or lies on some $\partial 3D_j$ that intersects $C(0,R-1)$.\\
		\textbf{Estimate on $W_2$.} Any point $x\in \overline{B(0,R-1)}\setminus \cup (3D_j)$ is either in $4D_j$ for some $j$
		or $x\in \overline{B(0,R-1)}\setminus \cup (4D_j)$. In the first case $u$ does not change sign in $B(x,1)$ and therefore
		$$|\nabla u(x)|\leq A|u(x)| \leq A e^N.$$
		In the second case $B(x,1) \subset B(0,R)\setminus \cup (3D_j)$ and $|\nabla u(x)| \leq A e^N$.
		Thus 
		$$ \max_{W_2}|u|, \max_{W_2}|\nabla u| \leq A e^N.$$
		Note also that 
		$$ W_2 \subset \overline{B(0,R-1)}\setminus B(0,R-7).$$

		Now, in $\Omega$, consider the analytic function 
		$$f(z)=\frac{u_x-iu_y}{z^k}, \quad |f(z)|=\frac{|\nabla u(z)|}{|z|^k}.$$
		Since $$\sup\limits_{B(0,R-R/64)\setminus \cup (3D_j)}|u|=1> \sup_{B(0,r/2)\setminus \cup 3D_j}|u|,$$  $$\max\limits_{W_1}  |u| \leq \left( \frac{r}{R} \right)^{2k} < \frac{1}{2},$$
		and since any point in $$\Omega_1=B(0,R-R/64)\setminus \cup (3D_j)$$ can be connected with $W_1$ by a curve of length at most $4R$ within $\Omega_1$, we must
		have $$\sup_\Omega|\nabla u| \geq \sup_{\Omega_1}|\nabla u| \geq \frac{1}{8R}$$ and
		$$ \sup_{\Omega_1} |f|\geq \left( R-R/64\right)^{-k} \sup_{\Omega_1}|\nabla u| \geq \frac{1}{8R} \left( R-R/64 \right)^{-k}.$$
		However 
		$$ \max_{W_1}|f|\leq \max_{W_1}|\nabla u |\left(\frac{2}{r} \right)^k \leq \left(\frac{r}{R} \right)^{2k}\left(\frac{2}{r} \right)^k= \left(\frac{2r}{R} \right)^{k} R^{-k} \leq $$ $$ \leq 2^{-k} R^{-k} < \frac{1}{8R} R^{-k} < \sup_{\Omega_1} |f|\leq \sup_{\Omega} |f| $$
		and
		$$ \max_{W_2}|f|\leq A e^N \frac{1}{(R-7)^k} = A e^N \left( \frac{R-R/64}{R-7}\right)^k \left( R-R/64\right)^{-k}\leq 
		$$ $$\leq A e^N \left( \frac{126}{127}\right)^k \left( R-R/64\right)^{-k} < \frac{1}{8R} \left( R-R/64\right)^{-k} \leq \sup_{\Omega} |f| $$
		if $R-7> \frac{127}{128}R $ and $C$ in the definition of $k$ is large enough. 
		By the maximum principle for holomorphic functions $\sup_\Omega |f|$ is achieved on $\partial 3D_j$ for some 
		$3D_j\subset B(0,R-1)\setminus \overline{B(0,r/2)}$. 
		
		For every disk $D_j$ with $3D_j\subset B(0,R-1)\setminus \overline{B(0,r/2)}$, consider the point $z_j$ on $\partial3D_j$ closest to the origin. All $3D_j$ that are not in the annulus  $B(0,R-1)\setminus \overline{B(0,r/2)}$ will not be considered further. Put $m_j=\min_{\partial3D_j}|u|$. Let $j_0$ be the index such that $$\frac{m_{j_0}}{|z_{j_0}|^k} = \max_j \frac{m_{j}}{|z_{j}|^k}.$$
		If $\sup_\Omega|f|$ is achieved on $\partial3D_{j_1}$, then for $x\in \Omega$,
		$$ \frac{|\nabla u(x)|}{|x|^k}\leq \frac{1}{|z_{j_1}|^k} \max_{\partial3D_{j_1}}|\nabla u|\leq A \frac{m_{j_1}}{|z_{j_1}|^k} \leq A \frac{m_{j_0}}{|z_{j_0}|^k}.$$
		So we conclude that 
		$$|\nabla u(x)| \leq \left( \frac{|x|}{|z_{j_0}|} \right)^k A m_{j_0}.$$
		Recalling that $\frac{1}{8R} \leq \sup_\Omega|\nabla u|$, we get
		$$ \frac{1}{8R} \leq \left(\frac{2R}{r}\right)^k Am_{j_0},$$
		so 
		$$ |u(z_{j_0})| \geq m_{j_0} \geq \frac{1}{8AR} \left(\frac{r}{2R}\right)^k.$$
		
		Now, let $(sz_{j_0},z_{j_0})$ be the longest subinterval of the radius $[0,z_{j_0}]$ starting at $z_{j_0}$ that is contained in $\Omega$. We have 
		$$|u(sz_{j_0})| \geq |u(z_{j_0})| - |z_{j_0}| \int_{s}^{1}|\nabla u(t z_{j_0})| dt \geq |u(z_{j_0})| - |z_{j_0}|   \int_{0}^{1}  Am_{j_0} t^k  dt \geq $$
		$$ \geq  m_{j_0}(1-\frac{AR}{k+1}) \geq \frac{m_{j_0}}{2}$$
		if the constant $C$ in the definition of $k$ is large enough.
		Note that $$\left(\frac{R}{r}\right)^{k} \geq 4^{k} > 16 AR \cdot 2^k.$$
		Hence $$\frac{m_{j_0}}{2} \geq \frac{1}{16AR} \left(\frac{r}{2R}\right)^k > \left(\frac{r}{R}\right)^{2k}$$ 
		and the point $sz_{j_0}$ cannot belong to $W_1$, whence it belongs to some $\partial 3D_j$ with $3D_j \subset B(0,R-1)\setminus \overline{B(0,r/2)}$.
		Then
		$$ \frac{m_j}{|z_j|^k} \geq \frac{|u(sz_{j_0})|}{A|sz_{j_0}|^k} \geq \frac{1}{2As^k} \frac{m_{j_0}}{|z_{j_0}|^k}.$$
		It remains to notice that, since the distance from $3D_j$ to $3D_{j_0}$ is at least $96$, we have 
		$$ s^k \leq (1-96/R)^k< \frac{1}{2A} $$
		if the constant $C$ in the definition of $k$ is large enough.
		But then $\frac{m_j}{|z_j|^k} > \frac{m_{j_0}}{|z_{j_0}|^k}$, which contradicts  the choice of $j_0$. Hence the hypothesis \ref{eq: wrong_hyp} is false.

	\end{proof}

	\section{Quasi-geodesic system of coordinate from \cite{Burq2022}}
	
\begin{lem}(quasi-geodesic coordinate)\\
	\label{lem: local_pseudogeod_coord_proof}
	Assume that $M$ is a $W^{2,\infty}$ Riemannian surface with boundary (the regularity here means that changes of charts are in $W^{2,\infty}\cap \mathcal{C}^1$), endowed with a Lipschitz (definite positive) metric $g$ and a Lipschitz (positive) density $K$. Let
	$$ \Delta = \frac{1}{K}\mathrm{div}(g^{-1}K\nabla).$$
	Then near any point of the boundary $x_0\in \partial M$ there exists a $W^{2,\infty}$ coordinate system $\Psi_0:U_0\subset \overline{M} \to B(0,\delta)\cap (\R_x \times[0, +\infty))$ with $\delta>0$ such that
	in this coordinate system $x_0 = (0, 0) $ and $$\Delta=\frac{1}{K}\nabla^{T}K\,b(x,y)\nabla,\quad b_{|\partial \Omega}=\begin{pmatrix}
		b_0(x) & 0\\
		0 & 1
	\end{pmatrix}$$
	with $b_0\in W^{1,\infty}(\partial\Omega)$ positive. 
\end{lem}
	
	\begin{proof}
		The proof is taken from the proposition $3.4$ of \cite{Burq2022}.
		
		Consider a point $x_0\in \partial M$. There exists a covering $ \partial M\subset \bigcup U_j$ , (here $\partial M$ is seen as a subset of $\overline{M}$ ) where $U_j$ are open sets of $\overline{M}$ and $W^{2,\infty}$ diffeomorphisms $$\Psi_j: V_j \to \R^2=\R_x\times \R_y,\quad  \overline{U_j} \subset V_j \subset M$$ such that $\Psi_{j}(V_{j})= B(0,1)\cap \{ y\geq 0 \}$ and $\Psi_j(U_j)\subset [-\delta,\delta]\times [0,\epsilon]$  for some $\delta,\, \epsilon >0$ small enough.
		Let $a=a(x,y)$ be the metric in this coordinate system, which is hence $W^{1,\infty}$ and defined for $||x|| \leq \delta',\, y\in [0,\epsilon']$, with $\delta'<\delta<1$ and $\epsilon'<\epsilon$. For any $x \in \{y=0\}$, consider the vector defined by
		\begin{equation}
			\label{eq: normalvectortoboundary}
			n(x) = (\lambda(x))^{-1/2}a^{-1}(0, x) \begin{pmatrix}
				1 \\
				0 
			\end{pmatrix},\quad \mbox{ for } \lambda=(1,0)\cdot a^{-1}(0,x) \begin{pmatrix}
				1 \\
				0 
			\end{pmatrix}.
		\end{equation}
		Because $a$ is $W^{1,\infty}$ and invertible, it is clear that  $n \in W^{1,\infty}([-\delta',\delta'])$ is the inward normal to the boundary for the metric a at the point $(0,x) \in \Psi_j (\partial M \cap V_j)$. Indeed, $$n(x)^{T} a(0, x) n(x) = (\lambda(x))^{-1}(1, 0)a^{-1}(0, x)a(0, x)a^{-1}(0, x)\begin{pmatrix}
			1 \\
			0 
		\end{pmatrix}= 1,$$
		which makes $n$ unitary and if $X\in \R_x$, then
		$$n(x)^{T}a(0, x)\begin{pmatrix}
			0 \\
			X 
		\end{pmatrix}= (\lambda(x))^{-1/2}(1, 0)a^{-1}(0, x)a(0, x)\begin{pmatrix}
			0 \\
			X 
		\end{pmatrix}=0,$$
		which proves that $n$ is orthogonal to the vectors tangent to the boundary. Finally, since its first component is positive, $n$ points inward.
		We now study the regularity of the quasi-geodesic coordinates. Let be $\chi \in \mathcal{C}_0^\infty([-\delta',\delta'])$ be equal to $1$ on $[-\delta'',\delta'']$ ($\delta''< \delta'$), and 
		$$m(x,y)=e^{-\langle y \mathcal{D}_x \rangle+1}(\chi n)(x,0),\quad \mbox{ with } \langle y\mathcal{D}_x \rangle= (1+|y\mathcal{D}_x|^2 )^{1/2}.$$ 
		
			\begin{lem}
			\label{lem: regularization_protocol}
			For any $q\in \R$ and $p>0$, the operators
			$$\langle y \mathcal{D}_x \rangle^{q}e^{-\langle y \mathcal{D}_x \rangle+1},\quad (y\mathcal{D}_{x})^p \langle y \mathcal{D}_x \rangle^{q}e^{-\langle y \mathcal{D}_x \rangle+1}$$
			are uniformly bounded on $L^{\infty}(\R,\R)$ with respect to $y$.
		\end{lem}
		\begin{proof}
			This come from the fact that this operator are convolution type operators with kernels
			$$K_{1,q}(z)=\frac{1}{y}\mathcal{F}\left(\langle \xi\rangle^q e^{-\langle \xi\rangle+1}\right)\left(\frac{z}{y}\right)$$
			and
			$$K_{2,q}(z)=\frac{1}{y}\mathcal{F}\left(\xi^p\langle\xi\rangle^q e^{-\langle \xi\rangle+1}\right)\left(\frac{z}{y}\right),$$
			where $\mathcal{F}$ is the Fourier transform on $\R$. The function $\xi \rightarrow \langle \xi\rangle^q e^{-\langle \xi\rangle+1}$ and $\xi \rightarrow i\xi\langle \xi\rangle^q e^{-\langle \xi\rangle+1}$ are Schwartz functions, Thus the Kernels are in $L^{1}(\R)$ with respect to $y$. Consequently, the operators are continuous on $L^{\infty}(\R)$.
		\end{proof}

		Because $n\in W^{1,\infty}(-\delta',\delta')$ and $e^{-\langle y \mathcal{D}_x \rangle+1}$ is uniformly bounded (with respect to $y$) on $L^\infty(\R_x)$, the map $(x, y) \to m(x, y)$ is Lipschitz and therefore the map $x \to m(x,0)$ is
		also Lipschitz. Moreover the map 
		$$
		(x, y) \in [-\epsilon', \epsilon'] \times [-\delta',\delta'] \to \phi_j (x, y) = x + y\,m(x, y) = x + y\,e^{-\langle y\mathcal{D}_x \rangle+1}(\chi\, n_{|s=0})$$
		is $W^{2,\infty}$. Indeed since  $(\chi\,n) \in L^\infty$ and $\partial_x (\chi\,n) \in L^\infty$ and thanks to the lemma \ref{lem: regularization_protocol}:
		
		$$ \partial_x \phi_j (x, y)=1+y\,e^{-\langle y\mathcal{D}_x \rangle+1}\partial_x(\chi\, n_{|s=0})\in L^\infty([-\epsilon', \epsilon'] \times [-\delta',\delta']),$$
		
		$$\partial_y \phi_j (x, y)=e^{-\langle y\mathcal{D}_x \rangle+1}(\chi\, n |s=0)-y^2\mathcal{D}_x^2\langle y\mathcal{D}_x\rangle^{-1} e^{-\langle y\mathcal{D}_x \rangle+1}(\chi\, n |s=0)$$
		$$=(1-\langle y\mathcal{D}_x \rangle+\langle y\mathcal{D}_x \rangle^{-1})e^{-\langle y\mathcal{D}_x \rangle+1}(\chi\, n |s=0)\in L^{\infty}([-\epsilon', \epsilon'] \times [-\delta',\delta']),$$
		
		$$\partial_x^2 \phi_j (x, y)=y\partial_x\,e^{-\langle y\mathcal{D}_x \rangle+1}\partial_x(\chi\, n_{|y=0})\in L^\infty([-\epsilon', \epsilon'] \times [-\delta',\delta']),$$ 
		
		$$\partial_y\partial_x \phi_j (x, y)=(1-\langle y\mathcal{D}_x \rangle+\langle y\mathcal{D}_x \rangle^{-1})\,e^{-\langle y\mathcal{D}_x \rangle+1}\partial_x(\chi\, n_{|s=0})\in L^\infty([-\epsilon', \epsilon'] \times [-\delta',\delta']),$$
		
		and
		$$\partial_y^2 \phi_j (x, y)=y\mathcal{D}_x(1-2\langle y\mathcal{D}_x\rangle^{-1}-\langle y \mathcal{D}_x\rangle^{-2}-\langle y \mathcal{D}_x\rangle^{-3})e^{-\langle y\mathcal{D}_x \rangle+1}\mathcal{D}_x(\chi\, n_{|s=0})\in L^\infty([-\epsilon', \epsilon'] \times [-\delta',\delta']).$$

		The differential of $\phi_j$ at $y = 0$ is
		\begin{equation}
			d_{x,y}{\phi_j}_{y=0}=
			\begin{pmatrix}
				 1 & \chi n(x, 0)_x \\
				 0 & \chi n(x,0)_y
			\end{pmatrix},
		\end{equation}
		which, according to \ref{eq: normalvectortoboundary} and the fact that $a$ is definite positive, is invertible for $x\in ]-\delta',\delta'[$.
		Hence, we deduce that $\phi_j$ is a $W^{2,\infty}$ diffeomorphism from a neighbourhood of $[-\delta',\delta']\times {0}_y$ to a
		neighbourhood of $[-\delta',\delta']\times {0}_y$ . Notice also that since $\phi_j$ sends the half plane ${y > 0}$ to itself, its
		inverse also sends the half plane ${s > 0}$ to itself. As a consequence, shrinking $U_j$ into a possibly smaller $U'_j$, we get a covering
		$$\partial M \subset \bigcup_{j=1}^{N}U'_j$$
		and $W^{2,\infty}$ diffeomorphisms $\Psi'_j=\phi_j^{-1}\circ \Psi_j$,
		such that after this change of variable, the metric $b(x, y)$
		is given for $y\geq0$ by
		$$b(x, y) ={\Psi'_j}^T a(x,y)\Psi'_j.$$
		In particular, for $y = 0^{+}$ we get
		\begin{equation}
			b(x,0^{+})=
			\begin{pmatrix}
				1 & 0 \\
				n(x)_x & n(x)_y 
			\end{pmatrix}
			a(x,0)
			\begin{pmatrix}
				1 &   n(x)_x \\
				0 &  n(x)_y
			\end{pmatrix}.
		\end{equation}
		Since $n(x)$ is the normal to the boundary we have
		$$n(x)^Ta(x, 0)n(x) = 1, (X, 0) a(x, 0)n(x) = 0, \forall\,X \in \R_x.$$
		We deduce
		\begin{equation}
			b(x,0^{+})=
			\begin{pmatrix}
				b_0(x) &   0\\
				0 &  1
			\end{pmatrix}
		\end{equation}
		with $b_0(x)>0$ Lipschitz.
		\end{proof}
	
	\bibliographystyle{elsarticle-harv} % aa 
	\typeout{}
	\bibliography{Source_Robin_Fort}

\begin{thebibliography}{21}
\expandafter\ifx\csname natexlab\endcsname\relax\def\natexlab#1{#1}\fi
\providecommand{\url}[1]{\texttt{#1}}
\providecommand{\href}[2]{#2}
\providecommand{\path}[1]{#1}
\providecommand{\DOIprefix}{doi:}
\providecommand{\ArXivprefix}{arXiv:}
\providecommand{\URLprefix}{URL: }
\providecommand{\Pubmedprefix}{pmid:}
\providecommand{\doi}[1]{\href{http://dx.doi.org/#1}{\path{#1}}}
\providecommand{\Pubmed}[1]{\href{pmid:#1}{\path{#1}}}
\providecommand{\bibinfo}[2]{#2}
\ifx\xfnm\relax \def\xfnm[#1]{\unskip,\space#1}\fi
%Type = Book
\bibitem[{Ahlfors(2006)}]{ahlfors2006lectures}
\bibinfo{author}{Ahlfors, L.V.}, \bibinfo{year}{2006}.
\newblock \bibinfo{title}{Lectures on quasiconformal mappings}.
  volume~\bibinfo{volume}{38}.
\newblock \bibinfo{publisher}{American Mathematical Soc.}
%Type = Book
\bibitem[{Astala et~al.(2008)Astala, Iwaniec and Martin}]{Astala2008-mi}
\bibinfo{author}{Astala, K.}, \bibinfo{author}{Iwaniec, T.},
  \bibinfo{author}{Martin, G.}, \bibinfo{year}{2008}.
\newblock \bibinfo{title}{Elliptic partial differential equations and
  quasiconformal mappings in the plane ({PMS-48})}.
\newblock Princeton Mathematical Series, \bibinfo{publisher}{Princeton
  University Press}, \bibinfo{address}{Princeton, NJ}.
%Type = Article
\bibitem[{Bourgain and Kenig(2005)}]{Bourgain2005}
\bibinfo{author}{Bourgain, J.}, \bibinfo{author}{Kenig, C.E.},
  \bibinfo{year}{2005}.
\newblock \bibinfo{title}{On localization in the continuous anderson-bernoulli
  model in higher dimension}.
\newblock \bibinfo{journal}{Inventiones mathematicae} \bibinfo{volume}{161},
  \bibinfo{pages}{389–426}.
\newblock \URLprefix \url{http://dx.doi.org/10.1007/s00222-004-0435-7},
  \DOIprefix\doi{10.1007/s00222-004-0435-7}.
%Type = Article
\bibitem[{Burq and Moyano(2022)}]{Burq2022}
\bibinfo{author}{Burq, N.}, \bibinfo{author}{Moyano, I.}, \bibinfo{year}{2022}.
\newblock \bibinfo{title}{Propagation of smallness and control for heat
  equations}.
\newblock \bibinfo{journal}{Journal of the European Mathematical Society}
  \bibinfo{volume}{25}, \bibinfo{pages}{1349–1377}.
\newblock \URLprefix \url{http://dx.doi.org/10.4171/JEMS/1213},
  \DOIprefix\doi{10.4171/jems/1213}.
%Type = Misc
\bibitem[{Das and Pinchover(2024)}]{das2024landisconjectureliouvillecomparison}
\bibinfo{author}{Das, U.}, \bibinfo{author}{Pinchover, Y.},
  \bibinfo{year}{2024}.
\newblock \bibinfo{title}{The landis conjecture via liouville comparison
  principle and criticality theory}.
\newblock \URLprefix \url{https://arxiv.org/abs/2405.11695},
  \href{http://arxiv.org/abs/2405.11695}{{\tt arXiv:2405.11695}}.
%Type = Misc
\bibitem[{Datchev(2014)}]{datchev2014quantitativelimitingabsorptionprinciple}
\bibinfo{author}{Datchev, K.}, \bibinfo{year}{2014}.
\newblock \bibinfo{title}{Quantitative limiting absorption principle in the
  semiclassical limit}.
\newblock \URLprefix \url{http://dx.doi.org/10.1007/s00039-014-0273-8},
  \DOIprefix\doi{10.1007/s00039-014-0273-8}.
%Type = Article
\bibitem[{Davey et~al.(2017)Davey, Kenig and Wang}]{Davey2017}
\bibinfo{author}{Davey, B.}, \bibinfo{author}{Kenig, C.},
  \bibinfo{author}{Wang, J.N.}, \bibinfo{year}{2017}.
\newblock \bibinfo{title}{The landis conjecture for variable coefficient
  second-order elliptic pdes}.
\newblock \bibinfo{journal}{Transactions of the American Mathematical Society}
  \bibinfo{volume}{369}, \bibinfo{pages}{8209–8237}.
\newblock \URLprefix \url{http://dx.doi.org/10.1090/tran/7073},
  \DOIprefix\doi{10.1090/tran/7073}.
%Type = Book
\bibitem[{Dyatlov and Zworski(2019)}]{dyatlovMathematicalTheoryScattering2019}
\bibinfo{author}{Dyatlov, S.}, \bibinfo{author}{Zworski, M.},
  \bibinfo{year}{2019}.
\newblock \bibinfo{title}{Mathematical {{Theory}} of {{Scattering
  Resonances}}}. volume \bibinfo{volume}{200} of
  \textit{\bibinfo{series}{Graduate {{Studies}} in {{Mathematics}}}}.
\newblock \bibinfo{publisher}{American Mathematical Society},
  \bibinfo{address}{Providence, Rhode Island}.
\newblock \DOIprefix\doi{10.1090/gsm/200}.
%Type = Misc
\bibitem[{{Fern{\'a}ndez-Bertolin} et~al.(2024){Fern{\'a}ndez-Bertolin}, Stan
  and Roncal}]{fernandez-bertolinLandisConjectureSurvey2024}
\bibinfo{author}{{Fern{\'a}ndez-Bertolin}, A.}, \bibinfo{author}{Stan, D.},
  \bibinfo{author}{Roncal, L.}, \bibinfo{year}{2024}.
\newblock \bibinfo{title}{Landis' conjecture: A survey}.
\newblock \DOIprefix\doi{10.48550/arXiv.2412.00788},
  \href{http://arxiv.org/abs/2412.00788}{{\tt arXiv:2412.00788}}.
%Type = Misc
\bibitem[{Frank and
  Ivanisvili(2026)}]{frank2026counterexampleslandisconjecturedimensions}
\bibinfo{author}{Frank, R.L.}, \bibinfo{author}{Ivanisvili, P.},
  \bibinfo{year}{2026}.
\newblock \bibinfo{title}{Counterexamples to the landis conjecture in
  dimensions three and higher}.
\newblock \URLprefix \url{https://arxiv.org/abs/2608.00802},
  \href{http://arxiv.org/abs/2608.00802}{{\tt arXiv:2608.00802}}.
%Type = Book
\bibitem[{Gilbarg and Trudinger(2001)}]{Gilbarg2001}
\bibinfo{author}{Gilbarg, D.}, \bibinfo{author}{Trudinger, N.S.},
  \bibinfo{year}{2001}.
\newblock \bibinfo{title}{Elliptic Partial Differential Equations of Second
  Order}.
\newblock \bibinfo{publisher}{Springer Berlin Heidelberg}.
\newblock \URLprefix \url{http://dx.doi.org/10.1007/978-3-642-61798-0},
  \DOIprefix\doi{10.1007/978-3-642-61798-0}.
%Type = Book
\bibitem[{Helffer(1988)}]{Helffer1988}
\bibinfo{author}{Helffer, B.}, \bibinfo{year}{1988}.
\newblock \bibinfo{title}{Semi-Classical Analysis for the Schr\"{o}dinger
  Operator and Applications}.
\newblock \bibinfo{publisher}{Springer Berlin Heidelberg}.
\newblock \URLprefix \url{http://dx.doi.org/10.1007/BFb0078115},
  \DOIprefix\doi{10.1007/bfb0078115}.
%Type = Article
\bibitem[{Klopp and Vogel(2019)}]{kloppSemiclassicalResolventEstimates2019a}
\bibinfo{author}{Klopp, F.}, \bibinfo{author}{Vogel, M.}, \bibinfo{year}{2019}.
\newblock \bibinfo{title}{Semiclassical resolvent estimates for bounded
  potentials}.
\newblock \bibinfo{journal}{Pure and Applied Analysis} \bibinfo{volume}{1},
  \bibinfo{pages}{1--25}.
\newblock \DOIprefix\doi{10.2140/paa.2019.1.1}.
%Type = Inbook
\bibitem[{Kondrat'ev and Landis(1991)}]{Kondratev1991}
\bibinfo{author}{Kondrat'ev, V.A.}, \bibinfo{author}{Landis, E.M.},
  \bibinfo{year}{1991}.
\newblock \bibinfo{title}{Qualitative Theory of Second Order Linear Partial
  Differential Equations}. \bibinfo{publisher}{Springer-Verlag},
  \bibinfo{address}{Berlin}. volume~\bibinfo{volume}{32} of
  \textit{\bibinfo{series}{Partial Differential Equations III}}.
\newblock pp. \bibinfo{pages}{87--192}.
%Type = Article
\bibitem[{Laurent and
  L{\'e}autaud(2023)}]{laurentUniformObservationSemiclassicalSchrodinger2023}
\bibinfo{author}{Laurent, C.}, \bibinfo{author}{L{\'e}autaud, M.},
  \bibinfo{year}{2023}.
\newblock \bibinfo{title}{Uniform observation of {{semiclassical
  Schr{\"o}dinger}} eigenfunctions on an interval}.
\newblock \bibinfo{journal}{Tunisian Journal of Mathematics}
  \bibinfo{volume}{5}, \bibinfo{pages}{125--170}.
\newblock \DOIprefix\doi{10.2140/tunis.2023.5.125}.
%Type = Article
\bibitem[{Logunov et~al.(2025)Logunov, Malinnikova, Nadirashvili and
  Nazarov}]{logunov2020landisconjectureexponentialdecay}
\bibinfo{author}{Logunov, A.}, \bibinfo{author}{Malinnikova, E.},
  \bibinfo{author}{Nadirashvili, N.}, \bibinfo{author}{Nazarov, F.},
  \bibinfo{year}{2025}.
\newblock \bibinfo{title}{The landis conjecture on exponential decay}.
\newblock \bibinfo{journal}{Inventiones mathematicae} \bibinfo{volume}{241},
  \bibinfo{pages}{465–508}.
\newblock \URLprefix \url{http://dx.doi.org/10.1007/s00222-025-01340-1},
  \DOIprefix\doi{10.1007/s00222-025-01340-1}.
%Type = Article
\bibitem[{Meshkov(1992)}]{Meshkov1992}
\bibinfo{author}{Meshkov, V.Z.}, \bibinfo{year}{1992}.
\newblock \bibinfo{title}{On the possible rate of decay at infinity of
  solutions of second order partial differential equations}.
\newblock \bibinfo{journal}{Mathematics of the USSR-Sbornik}
  \bibinfo{volume}{72}, \bibinfo{pages}{343–361}.
\newblock \URLprefix \url{http://dx.doi.org/10.1070/SM1992v072n02ABEH001414},
  \DOIprefix\doi{10.1070/sm1992v072n02abeh001414}.
%Type = Article
\bibitem[{Michel and Zworski(2018)}]{Kramers-Smoluchowski}
\bibinfo{author}{Michel, L.}, \bibinfo{author}{Zworski, M.},
  \bibinfo{year}{2018}.
\newblock \bibinfo{title}{A semiclassical approach to the kramers--smoluchowski
  equation}.
\newblock \bibinfo{journal}{SIAM Journal on Mathematical Analysis}
  \bibinfo{volume}{50}, \bibinfo{pages}{5362--5379}.
\newblock \URLprefix \url{https://doi.org/10.1137/17M1124826},
  \DOIprefix\doi{10.1137/17M1124826}.
%Type = Misc
\bibitem[{Shapiro(2019)}]{shapiro2017semiclassicalresolventboundsdimension}
\bibinfo{author}{Shapiro, J.}, \bibinfo{year}{2019}.
\newblock \bibinfo{title}{Semiclassical resolvent bounds in dimension two}.
\newblock \URLprefix \url{https://www.jstor.org/stable/26677715}.
%Type = Article
\bibitem[{Vodev(2020)}]{vodevSemiclassicalResolventEstimates2020}
\bibinfo{author}{Vodev, G.}, \bibinfo{year}{2020}.
\newblock \bibinfo{title}{Semiclassical resolvent estimates for
  {{H{\"o}lderpotentials}}}.
\newblock \bibinfo{journal}{Pure and Applied Analysis} \bibinfo{volume}{2},
  \bibinfo{pages}{841--860}.
\newblock \DOIprefix\doi{10.2140/paa.2020.2.841}.
%Type = Article
\bibitem[{Witten(1982)}]{10.4310/jdg/1214437492}
\bibinfo{author}{Witten, E.}, \bibinfo{year}{1982}.
\newblock \bibinfo{title}{{Supersymmetry and Morse theory}}.
\newblock \bibinfo{journal}{Journal of Differential Geometry}
  \bibinfo{volume}{17}, \bibinfo{pages}{661 -- 692}.
\newblock \URLprefix \url{https://doi.org/10.4310/jdg/1214437492},
  \DOIprefix\doi{10.4310/jdg/1214437492}.

\end{thebibliography}
	
\end{document}